\documentclass[leqno,a4paper,10pt]{amsart}
\usepackage{amsmath}
\usepackage{epsfig,graphics}
\usepackage{amssymb,latexsym,amscd}
\usepackage{mathrsfs}
\newtheorem{thm}{\bf Theorem}[section]
\newtheorem{df}[thm]{\bf Definition}
\newtheorem{prop}[thm]{\bf Proposition}
\newtheorem{cor}[thm]{\bf Corollary}
\newtheorem{lem}[thm]{\bf Lemma}
\newtheorem{rem}[thm]{\bf Remark}
\newtheorem{ex}[thm]{\bf Example}
\newcommand{\bs}{\boldsymbol}

\newcommand{\pf}{\noindent{\bfseries Proof. }}
\newcommand{\bi}{\bs{\rm i}}
\newcommand{\bj}{\bs{\rm j}}

\newcommand{\bk}{\bs{\rm k}}
\newcommand{\bl}{\bs{\rm l}}

\newcommand{\re}{{e}}
\newcommand{\rf}{{f}}
\newcommand{\B}{{\bf B}}
\newcommand{\M}{{\bf M}}

\newcommand{\ov}{\overline}
\newcommand{\w}{{\bf w}}

\numberwithin{equation}{section}

\begin{document}
\title[Crystal graphs and Cauchy decomposition]
{Crystal graphs for Lie superalgebras and Cauchy decomposition}
\author{JAE-HOON KWON}
\address{Department of Mathematics \\ University of Seoul \\ 90
Cheonnong-dong, Dongdaemun-gu \\ Seoul 130-743, Korea }
\email{jhkwon@uos.ac.kr }

\thanks{This research was supported by 2005 research fund of University of Seoul.} 

\maketitle

\begin{abstract}
We discuss Cauchy type decompositions of crystal graphs for general
linear Lie superalgebras. More precisely, we consider bicrystal
graph structures on various sets of matrices of non-negative
integers, and obtain their decompositions with explicit
combinatorial isomorphisms.
\end{abstract}

\section{introduction}
Let $\frak{gl}_{m|n}$ be the general linear Lie superalgebra over
$\mathbb{C}$, and let $\mathbb{C}^{m|n}$ be its $(m+n)$-dimensional
natural representation  with the $\mathbb{Z}_2$-grading
$(\mathbb{C}^{m|n})_0=\mathbb{C}^m$ and
$(\mathbb{C}^{m|n})_1=\mathbb{C}^n$. A tensor power of
$\mathbb{C}^{m|n}$ is completely reducible from Schur-Weyl duality,
and its irreducible components, called irreducible polynomial
representations, are parameterized by $\mathcal{P}_{m|n}$, the set
of all $(m,n)$-hook partitions \cite{BR}. Let
$S(\mathbb{C}^{m|n}\otimes\mathbb{C}^{u|v})$ be the (super)
symmetric algebra generated by the
$(\frak{gl}_{m|n},\frak{gl}_{u|v})$-bimodule
$\mathbb{C}^{m|n}\otimes\mathbb{C}^{u|v}$. From Howe duality, it is
also completely reducible as a
$(\frak{gl}_{m|n},\frak{gl}_{u|v})$-bimodule, and we have the
following Cauchy type decomposition;
\begin{equation}\label{Cauchydecomp}
S(\mathbb{C}^{m|n}\otimes\mathbb{C}^{u|v})=\bigoplus_{\lambda\in\mathcal{P}_{m|n}\cap
\mathcal{P}_{u|v}} V_{m|n}(\lambda)\otimes V_{u|v}(\lambda),
\end{equation}
where $V_{m|n}(\lambda)$ and $V_{u|v}(\lambda)$ denote the
irreducible polynomial representations of $\frak{gl}_{m|n}$ and
$\frak{gl}_{u|v}$, respectively, corresponding to $\lambda$ (see
\cite{CW,H}). In terms of characters, \eqref{Cauchydecomp} also
yields a Cauchy type identity of hook Schur functions
(cf.\cite{Mac}).

The purpose of this paper is to understand the decomposition
\eqref{Cauchydecomp} within a framework of (abstract) crystal graphs
for Lie superalgebras which were developed by Benkart, Kang and
Kashiwara \cite{BKK}. For $\lambda\in\mathcal{P}_{m|n}$, we denote
by $\B_{m|n}(\lambda)$ the set of all $(m,n)$-hook semistandard
tableaux of shape $\lambda$, which parameterizes the basis element
of $V_{m|n}(\lambda)$ \cite{BR}. According to the crystal base
theory in \cite{BKK}, $\B_{m|n}(\lambda)$ becomes a colored oriented
graph, which we call a {\it crystal graph for $\frak{gl}_{m|n}$} or
{\it $\frak{gl}_{m|n}$-crystal}. As in the case of symmetrizable
Kac-Moody algebras, the crystal graphs for $\frak{gl}_{m|n}$ have
nice behaviors under tensor product, and we can decompose various
finite dimensional representations of $\frak{gl}_{m|n}$ in a purely
combinatorial way (cf.\cite{KK}).

For non-negative integers $m,n,u,v$ such that $m+n,u+v>0$, let
\begin{equation*}
\M=\{\,A=(a_{bb'})_{b\in \B_{m|n}, b'\in \B_{u|v}}\,|\, \text{(i)
$a_{bb'}\in\mathbb{Z}_{\geq 0}$, (ii) $a_{bb'}\leq 1$ if
$|b|\neq|b'|$}\, \},
\end{equation*}
where $\B_{m|n}$ (resp.$\B_{u|v}$) is the crystal graph associated
to the natural representation of $\frak{gl}_{m|n}$ (resp.
$\frak{gl}_{u|v}$), and $|b|$ denotes the degree of $b$. Note that
$\M$ naturally parameterizes the set of monomial basis of
$S(\mathbb{C}^{m|n}\otimes\mathbb{C}^{u|v})$. Then we show that as a
crystal graph for $\frak{gl}_{m|n}\oplus \frak{gl}_{u|v}$ (or
$(\frak{gl}_{m|n},\frak{gl}_{u|v})$-bicrystal),
\begin{equation}\label{main result1}
\M\simeq \bigoplus_{\lambda\in
\mathcal{P}_{m|n}\cap\mathcal{P}_{u|v}} \B_{m|n}(\lambda)\times
\B_{u|v}(\lambda).
\end{equation}
The isomorphism is given by a super-analogue of the well-known Knuth
correspondence \cite{Kn}. But our proof is  different from the
original one  since the decomposition is given by characterizing all
the highest weight elements of the connected components from the
view point of crystal graphs. Furthermore, our approach enables us
to explain several variations of the Knuth correspondence
(cf.\cite{Fu}) in a unified way, and also derive  an interesting
relation between the statistics of the diagonal entries of a
symmetric matrix in $\M$ and the number of odd parts in the shape of
the corresponding tableau (a special case of this relation was first
observed in \cite{Kn}).

We may naturally extend the above decomposition to a semi-infinite
case. Let $\frak{g}$ be a contragredient Lie superalgebra of
infinite rank whose Dynkin diagram is given by \begin{center}
\setlength{\unitlength}{0.45cm}
\begin{picture}(15,2)(0,0)
\put(-.85,.5){\line(1,0){1}}
\put(0,0){\makebox(1,1){$\bigcirc$}}\put(.85,.5){\line(1,0){1}}
\put(5,0){\makebox(1,1){$\bigcirc$}}\put(3.85,.5){\line(1,0){1.2}}\put(5.85,.5){\line(1,0){1.2}}
\put(7,0){\makebox(1,1){$\bigotimes$}}\put(7.85,.5){\line(1,0){1.2}}
\put(9,0){\makebox(1,1){$\bigcirc$}}\put(9.85,.5){\line(1,0){1.2}}
\put(14,0){\makebox(1,1){$\bigcirc$}}\put(12.85,.5){\line(1,0){1.2}}\put(14.85,.5){\line(1,0){1.2}}
\put(2.5,.2){$\cdots$}\put(11.5,.2){$\cdots$}
\put(-2,.2){$\cdots$}\put(16.5,.2){$\cdots$}
\end{picture}\ \ \ \ \ \ \ \ \ \ \
\end{center} \vskip .5cm
(see \cite{Kac}). First, we introduce a $\frak{g}$-crystal
$\mathscr{F}$ consisting of semi-infinite words, which can be viewed
as a crystal graph associated to a Fock space representation of
$\frak{g}$ analogous to the level one fermionic Fock space
representation of $\widehat{\frak{gl}}_{\infty}$ (cf.\cite{Kac90}).
We show that a connected component of a tensor power
$\mathscr{F}^{\otimes u}$ ($u\geq 1$) can be realized as the set of
semi-infinite semistandard tableaux, which is generated by a highest
weight element. Moreover, an explicit multiplicity-free
decomposition of $\mathscr{F}^{\otimes u}$ as a
$(\frak{g},\frak{gl}_u)$-bicrystal is given, where each connected
component is parameterized by a generalized partition of length $u$.
To prove this, we identify $\mathscr{F}^{\otimes u}$ with a set of
certain matrices with infinite number of rows, and apply the methods
used in the case of finite ranks. More precisely, an element in
$\mathscr{F}^{\otimes u}$ is equivalent to a unique pair of an
semi-infinite semistandard tableau and a rational semistandard
tableau as an element in a $(\frak{g},\frak{gl}_u)$-bicrystal.
Hence, it gives rise to a Knuth correspondence of a semi-infinite
type, and a Cauchy type identity. As a by-product, we obtain a
character formula of a $\frak{g}$-crystal of semi-infinite
semistandard tableaux occurring as a connected component in
$\mathscr{F}^{\otimes u}$. This character formula is given in terms
of ordinary Schur functions and the Littlewood-Richardson
coefficients, and it is very similar to the ones of the irreducible
highest weight representations of $\widehat{\frak{gl}}_{\infty}$ or
$\widehat{\frak{gl}}_{\infty|\infty}$ obtained in \cite{CL,KacR}. In
fact, using crystal graphs, we can give a similar combinatorial
proof of the Cauchy type decomposition of a higher level Fock space
representation of $\widehat{\frak{gl}}_{\infty}$ given by Kac and
Radul \cite{KacR}, and hence the character formula of irreducible
highest weight representations of $\widehat{\frak{gl}}_{\infty}$. We
also expect a combinatorial proof of the decomposition of a Fock
space representation of $\widehat{\frak{gl}}_{\infty|\infty}$ given
by Cheng and Lam \cite{CL}.

The paper is organized as follows. In Section 2, we review the basic
notions and the main results on crystal graphs for $\frak{gl}_{m|n}$
in \cite{BKK}. In Section 3, we prove the decomposition \eqref{main
result1} and then study the diagonal action of the Kashiwara
operators on the set of symmetric matrices in $\M$. In Section 4, we
describe the dual decomposition which is associated to the (super)
exterior algebra $\Lambda(\mathbb{C}^{m|n}\otimes\mathbb{C}^{u|v})$.
Finally, in Section 5, we introduce a $\frak{g}$-crystal of
semi-infinite semistandard tableaux, and a $\frak{gl}_u$-crystal of
rational semistandard tableaux (cf.\cite{St}) . Then using these
combinatorial realizations of crystal graphs, we establish a Cauchy
type decomposition of a tensor power $\mathscr{F}^{\otimes u}$ as a
$(\frak{g},\frak{gl}_u)$-bicrystal. \vskip 3mm

{\bf Acknowledgment } The author would like to thank Prof. S.-J.
Kang and Prof. S.-J. Cheng for their interests in this work and many
helpful discussions.

\section{Crystal graphs for $\frak{gl}_{m|n}$}
In this section, we recall the basic notions on crystal graphs for
$\frak{gl}_{m|n}$ developed in \cite{BKK}.

\subsection{Definitions}

For non-negative integers $m,n$ with $m+n>0$, let $\frak{gl}_{m|n}$
be the general linear Lie superalgebra over $\mathbb{C}$ (see
\cite{Kac}). Let
$$\B_{m|n}=\{\, \ov{m} < \ov{m-1} < \cdots < \ov{1} < 1 < 2 <
\cdots < n\, \}$$ be a linearly ordered set. Set $\B_{m|n}^{+} =
\{\, \ov{m}, \ov{m-1}, \cdots , \ov{1}\, \}$ and $\B_{m|n}^{-} =
\{\, 1, 2,
 \cdots , n\, \}$. For $b\in
\B_{m|n}$, we define $|b|$, {\it degree of $b$}, by $|b|=0$ (resp.
$1$) if $b\in \B_{m|n}^{+}$ (resp. $\B_{m|n}^-$). The free abelian
group $P_{m|n}=\bigoplus_{b\in \B_{m|n}} \mathbb{Z} \epsilon_{b}$,
which is generated by $\epsilon_b$ ($b\in \B_{m|n}$), is called the
{\it weight lattice} of $\frak{gl}_{m|n}$. There is a natural
symmetric $\mathbb{Z}$-bilinear form $(\ ,\,)$ on $P_{m|n}$, where
$(\epsilon_b,\epsilon_{b'})=(-1)^{|b|}\delta_{bb'}$ for
$b,b'\in\B_{m|n}$. Let
\begin{equation*}
I_{m|n}=\{\,\ov{m-1},\cdots,\ov{1},0,1,\cdots,n-1 \,\}.
\end{equation*}
The {\it simple root} $\alpha_i$ ($i\in I_{m|n}$) of
$\frak{gl}_{m|n}$ is given by
\begin{equation*}
\begin{cases}
\alpha_{\overline{k}} =
\epsilon_{\overline{k+1}}-\epsilon_{\overline {k}}, & \text{for
$k=1,\cdots,m-1$}, \\
\alpha_l = \epsilon_l-\epsilon_{l+1}, & \text{for $l=1,\cdots,
n-1$},
\\ \alpha_0 = \epsilon_{\ov{1}}-\epsilon_{1}. &
\end{cases}
\end{equation*}
Set $Q=\bigoplus_{i\in I_{m|n}} \mathbb{Z} \alpha_i$, which we call
the {\it root lattice of $\frak{gl}_{m|n}$}. A partial ordering on
$P_{m|n}$ is given by $\lambda\geq \mu$ if and only if
$\lambda-\mu\in \sum_{i\in I_{m|n}} \mathbb{Z}_{\geq 0} \alpha_i$
for $\lambda,\mu\in P_{m|n}$. We also define the {\it simple coroot
$h_i\in P_{m|n}^*$} ($i\in I_{m|n}$) by
\begin{equation*}
\langle h_i,\lambda \rangle =
\begin{cases}
(\alpha_i,\lambda), & \text{if $i=\ov{m-1},\cdots,\ov{1},0$}, \\
-(\alpha_i,\lambda), & \text{if $i=1,\cdots,n-1$},
\end{cases}
\end{equation*}
for $\lambda\in P_{m|n}$, where $\langle\ ,\,\rangle$ is the natural
pairing on $P_{m|n}^*\times P_{m|n}$. With respect to the above
simple roots, the Dynkin diagram is
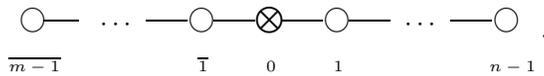
\begin{center}
\setlength{\unitlength}{0.45cm}
\begin{picture}(15,3)(0,0)
\put(0,0){\makebox(1,1){$\bigcirc$}}\put(.85,.5){\line(1,0){1}}
\put(5,0){\makebox(1,1){$\bigcirc$}}\put(3.85,.5){\line(1,0){1.2}}\put(5.85,.5){\line(1,0){1.2}}
\put(7,0){\makebox(1,1){$\bigotimes$}}\put(7.85,.5){\line(1,0){1.2}}
\put(9,0){\makebox(1,1){$\bigcirc$}}\put(9.85,.5){\line(1,0){1.2}}
\put(14,0){\makebox(1,1){$\bigcirc$}}\put(12.85,.5){\line(1,0){1.2}}
\put(2.5,.2){$\cdots$}\put(11.5,.2){$\cdots$}

\put(-.2,-1){\tiny $\ov{m-1}$}\put(5.4,-1){\tiny
$\ov{1}$}\put(7.4,-1){\tiny $0$}\put(9.4,-1){\tiny $1$}\put(14
,-1){\tiny $n-1$}
\end{picture}\ \ .
\end{center} \vskip .5cm

Motivated by the crystal bases of integral representations of the
quantized enveloping algebra $U_q(\frak{gl}_{m|n})$, we introduce
the notion of abstract crystal graphs for $\frak{gl}_{m|n}$.
\begin{df}\label{crystal graph}
{\rm (1) A  {\it crystal graph for $\frak{gl}_{m|n}$} (or {\it
$\frak{gl}_{m|n}$-crystal}) is a set $B$ together with the maps
\begin{equation*}
\begin{split}
&{\rm wt}  : B \rightarrow P_{m|n},  \\
&\varepsilon_i, \varphi_i: B \rightarrow
\mathbb{Z}_{\geq 0}, \\
&{e}_i, {f}_i: B \rightarrow
B\cup\{0\},
\end{split}
\end{equation*}
for $i\in I_{m|n}$ ($0$ is a formal symbol), satisfying the
following conditions;
\begin{itemize}

\item[(a)] for $i\in I_{m|n}$ and $b\in B$, we have
\begin{equation*}
\begin{aligned}\text{}
& \varphi_i(b)-  \varepsilon_i(b) =\langle h_i,{\rm wt}(b)\rangle, \ \ \ \text{($i\neq 0$)}, \\
& \langle h_0,{\rm wt}(b)\rangle\geq 0, \text{ and } \
\varphi_0(b)+ \varepsilon_0(b) =
\begin{cases}
0, & \text{if  $\langle h_0,{\rm wt}(b)\rangle=0$},\\
1, & \text{if  $\langle h_0,{\rm wt}(b)\rangle>0$},
\end{cases}
\end{aligned}
\end{equation*}

\item[(b)] if ${e}_i b \in B$ for $i\in I_{m|n}$ and $b\in B$, then
$$\varepsilon_i({e}_i b) = \varepsilon_i(b) - 1, \quad
\varphi_i({e}_i b) = \varphi_i(b) + 1, \quad {\rm wt}(e_ib)={\rm wt}(b)+\alpha_i,$$

\item[(c)] if ${f}_i b \in B$ for $i\in I_{m|n}$ and $b\in B$, then
$$\varepsilon_i({f}_i b) = \varepsilon_i(b) + 1, \quad
\varphi_i({f}_i b) = \varphi_i(b) - 1,\quad {\rm wt}(f_ib)={\rm wt}(b)-\alpha_i,$$

\item[(d)] ${f}_i b = b'$ if and only if $b = {e}_i
b'$ for all $i\in I_{m|n}$, $b, b' \in B$,
\end{itemize}
\ \ \ \ \ (We call $e_i$ and $f_i$ ($i\in I_{m|n}$) the {\it
Kashiwara operators}).\vskip 3mm

(2) Let $B$ be a crystal graph for $\frak{gl}_{m|n}$. A subset
$B'\subset B$ is called a {\it subcrystal of $B$} if $B'$ is itself
a crystal graph for $\frak{gl}_{m|n}$ with respect to ${\rm
wt},\varepsilon_i,\varphi_i,e_i,f_i$ ($i\in I_{m|n}$) of $B$.}
\end{df}\vskip 3mm
\begin{rem}{\rm
(1) The above definition is based on the crystal bases of integrable
representations of $U_q(\frak{gl}_{m|n})$ in \cite{BKK}, while
crystal graphs for contragredient Lie superalgebras might be defined
in a more general sense, as in the case of symmetrizable Kac-Moody
algebras (cf. \cite{Kas93,Kas94}).

(2) A crystal graph $B$ becomes an $I_{m|n}$-colored oriented graph,
where
\begin{equation*}
b\stackrel{i}{\rightarrow}b' \ \ \ \ \text{if and only if}  \ \ \ \
\ b'=f_ib \ \ \ \ \ \ (i\in I_{m|n}).
\end{equation*}
}
\end{rem}\vskip 3mm

\begin{df}\label{tensor product}{\rm
Let $B_1$ and $B_2$ be crystal graphs for $\frak{gl}_{m|n}$.
We define the {\it tensor product  of} $B_1$ and $B_2$ to be the set
$B_1\otimes B_2=\{\,b_1\otimes b_2\,|\,b_i\in B_i,\,\, (i=1,2)\,\}$
with
{\allowdisplaybreaks
\begin{equation*}
{\rm wt}(b_1\otimes b_2)={\rm wt}(b_1)+{\rm wt}(b_2),
\end{equation*}
\begin{equation*}
\varepsilon_i(b_1\otimes b_2)=
\begin{cases}
{\rm max}(\varepsilon_i(b_1),\varepsilon_i(b_2)-\langle h_i,{\rm
wt}(b_1)\rangle),
& \text{if $i=\ov{m-1},\cdots,\ov{1}$}, \\
{\rm max}(\varepsilon_i(b_1)-\langle h_i,{\rm
wt}(b_2)\rangle,\varepsilon_i(b_2)),
& \text{if $i=1,\cdots,n-1$}, \\
\varepsilon_0(b_1), & \text{if  $\langle h_0,{\rm wt}(b_1)\rangle>0$}, \\
\varepsilon_0(b_2), & \text{if  $\langle h_0,{\rm
wt}(b_1)\rangle=0$},
\end{cases}
\end{equation*}
\begin{equation*}
\varphi_i(b_1\otimes b_2)=
\begin{cases}
{\rm max}(\varphi_i(b_1)+\langle h_i,{\rm
wt}(b_2)\rangle,\varphi_i(b_2)),
& \text{if $i=\ov{m-1},\cdots,\ov{1}$}, \\
{\rm max}(\varphi_i(b_1),\varphi_i(b_2)+\langle h_i,{\rm
wt}(b_1)\rangle),
& \text{if $i=1,\cdots,n-1$}, \\
\varphi_0(b_1), & \text{if  $\langle h_0,{\rm wt}(b_1)\rangle>0$}, \\
\varphi_0(b_2), & \text{if  $\langle h_0,{\rm wt}(b_1)\rangle=0$},
\end{cases}
\end{equation*}
\begin{equation*}
e_i(b_1\otimes b_2)=
\begin{cases}
{e}_i b_1 \otimes b_2, & \text{if $i=\ov{m-1},\cdots,\ov{1}$,\ \
$\varphi_i(b_1)\geq \varepsilon_i(b_2)$}, \\
b_1\otimes {e}_i b_2, & \text{if $i=\ov{m-1},\cdots,\ov{1}$,\ \
$\varphi_i(b_1)< \varepsilon_i(b_2)$},\\
{e}_i b_1 \otimes b_2, & \text{if $i=1,\cdots,n-1$,\ \
$\varphi_i(b_2)< \varepsilon_i(b_1)$}, \\
b_1\otimes {e}_i b_2, & \text{if $i=1,\cdots,n-1$,\ \
$\varphi_i(b_2)\geq \varepsilon_i(b_1)$}, \\
{e}_0 b_1 \otimes b_2, & \text{if  $\langle h_0,{\rm wt}(b_1)\rangle>0$}, \\
 b_1 \otimes {e}_0 b_2, & \text{if  $\langle h_0,{\rm wt}(b_1)\rangle=0$},
\end{cases}
\end{equation*}
\begin{equation*}
f_i(b_1\otimes b_2)=
\begin{cases}
{f}_i b_1 \otimes b_2, & \text{if $i=\ov{m-1},\cdots,\ov{1}$,\ \
$\varphi_i(b_1)>\varepsilon_i(b_2)$}, \\
b_1\otimes {f}_i b_2, & \text{if $i=\ov{m-1},\cdots,\ov{1}$,\ \
$\varphi_i(b_1)\leq \varepsilon_i(b_2)$},\\
{f}_i b_1 \otimes b_2, & \text{if $i=1,\cdots,n-1$,\ \
$\varphi_i(b_2)\leq \varepsilon_i(b_1)$}, \\
b_1\otimes {f}_i b_2, & \text{if $i=1,\cdots,n-1$,\ \
$\varphi_i(b_2)> \varepsilon_i(b_1)$}, \\
{f}_0 b_1 \otimes b_2, & \text{if  $\langle h_0,{\rm wt}(b_1)\rangle>0$}, \\
 b_1 \otimes {f}_0 b_2, & \text{if  $\langle h_0,{\rm
 wt}(b_1)\rangle=0$},
\end{cases}
\end{equation*}}
where we assume that $0\otimes b_2=b_1\otimes 0=0$.}
\end{df}
Then, it is straightforward to check that $B_1\otimes B_2$ is a
crystal graph for $\frak{gl}_{m|n}$.

\begin{df}{\rm Let $B_1$ and $B_2$ be two crystal graphs for
$\frak{gl}_{m|n}$.
\begin{itemize}
\item[(1)] The {\it direct sum} $B_1\oplus B_2$ is the
disjoint union of $B_1$ and $B_2$.

\item[(2)] An {\it isomorphism} $\psi : B_1 \rightarrow B_2$ of
$\frak{gl}_{m|n}$-crystals is an isomorphism of $I_{m|n}$-colored
oriented graphs which preserves ${\rm wt}$, $\varepsilon_i$, and
$\varphi_i$ ($i\in I_{m|n}$). We say that {\it $B_1$ is isomorphic
to $B_2$}, and write $B_1\simeq B_2$.

\item[(3)] For $b_i\in B_i$ ($i=1,2$), let $C(b_i)$ denote the connected
component of $b_i$ as an $I_{m|n}$-colored oriented graph. We say
that $b_1$ {\it is $\frak{gl}_{m|n}$-equivalent to} $b_2$, if there
is an isomorphism of crystal graphs $C(b_1)\rightarrow C(b_2)$
sending $b_1$ to $b_2$. We often write $b_1
{\simeq}_{\frak{gl}_{m|n}} b_2$ (or simply $b_1 {\simeq} b_2$ if
there is no confusion).
\end{itemize}
}
\end{df}

\subsection{Semistandard tableaux for $\frak{gl}_{m|n}$}
$\B_{m|n}$ becomes a crystal graph for $\frak{gl}_{m|n}$ whose
associated $I_{m|n}$-colored oriented graph is given by
\begin{equation*}
\ov{m}\ \stackrel{\ov{m-1}}{\longrightarrow}\ \ov{m-1}\ \stackrel{\ov{m-2}}{\longrightarrow}
\cdots\stackrel{\ov{1}}{\longrightarrow}\ \ov{1}\ \stackrel{0}{\longrightarrow}\ 1 \
\stackrel{1}{\longrightarrow}\cdots\stackrel{n-2}{\longrightarrow} n-1\stackrel{n-1}{\longrightarrow} n,
\end{equation*}
where ${\rm wt}(b)=\epsilon_b$, and $\varepsilon_i(b)$ (resp.
$\varphi_i(b)$) is the number of $i$-colored arrows coming into $b$
(resp. going out of $b$) for $i\in I_{m|n}$ and $b\in\B_{m|n}$. Note
that $\B_{m|n}$ is the crystal graph associated to the natural
representation $\mathbb{C}^{m|n}$.

Let $\mathcal{W}_{m|n}$ be the set of all finite words with the
letters in $\B_{m|n}$. The empty word is denoted by $\emptyset$.
Then $\mathcal{W}_{m|n}$ is a crystal  graph for $\frak{gl}_{m|n}$
since we may identify each non-empty word $w=w_1\cdots w_r$ with
$w_1\otimes\cdots\otimes w_r\in \B_{m|n}^{\otimes r}$, where
$\{\emptyset\}$ forms a trivial crystal graph, that is, ${\rm
wt}(\emptyset)=0$, $e_i\emptyset=f_i\emptyset=0$, and
$\varepsilon_i(\emptyset)=\varphi_i(\emptyset)=0$ for all $i\in
I_{m|n}$.

Following the tensor product rule in Definition \ref{tensor
product}, we can describe the Kashiwara operators $\re_i, \rf_i :
\mathcal{W}_{m|n} \rightarrow \mathcal{W}_{m|n}\cup \{0\}$ ($i\in
I_{m|n}$) in a more explicit way;

\begin{itemize}
\item[(1)] Suppose that a non-empty word $w=w_1\cdots w_r$ is given.
To each letter $w_k$, we assign
\begin{equation*}
\epsilon^{(i)}(w_k)=
\begin{cases}
+, & \text{if ($i=\ov{p}$, $w_k=\ov{p+1}$), or ($i=p$, $w_k=p$), or ($i=0$, $w_k=\ov{1}$)}, \\
-, & \text{if ($i=\ov{p}$, $w_k=\ov{p}$), or ($i=p$, $w_k=p+1$), or ($i=0$, $w_k=1$)}, \\
\,\,\cdot\ , & \text{otherwise},
\end{cases}
\end{equation*}
and let
$\epsilon^{(i)}(w)=(\epsilon^{(i)}(w_1),\cdots,\epsilon^{(i)}(w_r))$.

\item[(2)] If $i=\ov{p}$ for $1\leq p\leq m-1$, then we replace a pair
$(\epsilon^{(i)}(w_s),\epsilon^{(i)}(w_{s'}))=(+,-)$ such that
$s<s'$ and $\epsilon^{(i)}(w_t)=\cdot$ for $s<t<s'$ by
$(\,\cdot\,,\,\cdot\,)$ in $\epsilon^{(i)}(w)$, and repeat this
process as far as possible until we get a sequence with no $+$
placed to the left of $-$. If $i=p$ for $1\leq p \leq n-1$, then we
do the same work for $(-,+)$-pair in $\epsilon^{(i)}(w)$ until we
get a sequence with no $-$ placed to the left of $+$. We call this
sequence the {\it $i$-signature} of $w$. If $i=0$, then we define
{\it $0$-signature of $w$} to be $\epsilon^{(i)}(w_k)(\neq \cdot)$
such that $\epsilon^{(i)}(w_l) = \cdot$ for all $1\leq l <k$.

\item[(3)] If $i=\ov{p}$ (resp. $i=p$), then we call the right-most
(resp. left-most) $-$ in the $i$-signature of $w$ the {\it $i$-good
$-$ sign}, and define $\re_i w$ to be the word obtained by applying
$e_i$ to $\ov{p}$ (resp. $p+1$) corresponding to the $i$-good $-$
sign. If there is no $i$-good $-$ sign, then we define $e_iw=0$.

\item[(4)] If $i=\ov{p}$ (resp. $i=p$), then we call the left-most
(resp. right-most) $+$ in the $i$-signature of $w$ the {\it $i$-good
$+$ sign}, and define $\rf_i w$ to be the word obtained by applying
$f_i$ to $\ov{p+1}$ (resp. $p$) corresponding to the $i$-good $+$
sign. If there is no $i$-good $+$ sign, then we define $f_iw=0$.

\item[(5)] If $i=0$, then we define $\re_0 w$ (resp. $\rf_0 w$) to be the word obtained
by applying $e_0$ (resp. $f_0$) to the letter corresponding to the
$0$-signature of $w$. If the $0$-signature of $w$ is empty, then we
define $\re_0 w= 0$ (resp. $f_0w=0$).
\end{itemize}
Note that we have
\begin{equation*}
\varepsilon_i(w)={\rm max}\{\,k\,|\,e_i^kw\neq 0\,\},
\ \ \ \varphi_i(w)={\rm max}\{\,k\,|\,f_i^kw\neq 0\,\},
\end{equation*}
for $w\in\mathcal{W}_{m|n}$ and $i\in I_{m|n}$. \vskip 3mm

\begin{ex}{\rm Suppose that
\begin{equation*}
w=\, 1 \ \ov{1} \ \ov{1} \ \ov{1} \ \ov{2} \ 2 \ 2 \ \ov{2} \ 1.
\end{equation*}
Then
\begin{equation*}
\begin{split}
&\epsilon^{(\ov{1})}(w)=(\ \cdot\ , \, - \, , \, - \, , \, \ominus \, , \,\oplus \,,\ \cdot \ ,\ \cdot \ ,\ + \,,\ \cdot \ ),  \\
&\epsilon^{(1)}(w)=(\, \oplus \,,\ \cdot \ ,\ \cdot \ ,\  \cdot\ , \ \cdot \ , \ \ominus \, , \, - \, , \ \cdot \ , \,  + \ ),  \\
&\epsilon^{(0)}(w)=(\, \ominus \, , \, + \,,\,+\,,\,+\,,\ \cdot \ ,\
\cdot \ ,\  \cdot\ , \ \cdot \ , \ - \ ),
\end{split}
\end{equation*}
where $\oplus,\ominus$ denote the $i$-good signs ($i\neq 0$), or the
$0$-signature. We have
\begin{equation*}
\begin{split}
&\re_{\ov{1}}(w)= \, 1 \ \ov{1} \ \ov{1} \ \ov{2} \ \ov{2} \ 2 \ 2 \ \ov{2} \ 1, \ \ \
\rf_{\ov{1}}(w)= \, 1 \ \ov{1} \ \ov{1} \ \ov{1} \ \ov{1} \ 2 \ 2 \ \ov{2} \ 1, \\
&\re_1(w)= \, 1 \ \ov{1} \ \ov{1} \ \ov{1} \ \ov{2} \ 1 \ 2 \ \ov{2} \ 1,\ \ \
\rf_1(w)= \, 2 \ \ov{1} \ \ov{1} \ \ov{1} \ \ov{2} \ 2 \ 2 \ \ov{2} \ 1, \\
&\re_0(w)= \, \ov{1} \ \ov{1} \ \ov{1} \ \ov{1} \ \ov{2} \ 2 \ 2 \ \ov{2} \ 1, \ \ \
\rf_0(w)= 0 .
\end{split}
\end{equation*}
}
\end{ex}\vskip 3mm

A {\it partition} is a non-increasing sequence of non-negative
integers $\lambda = (\lambda_k)_{k\geq 1}$ such that all but a
finite number of its terms are zero. Each $\lambda_k$ is called a
{\it part of $\lambda$}, and the number of non-zero parts is called
the {\it length of $\lambda$}. We also write
$\lambda=(1^{m_1},2^{m_2},\cdots)$ where $m_i$ is the number of
occurrences of $i$ in $\lambda$. Recall that a partition $\lambda =
(\lambda_k)_{k\geq 1}$ is identified with a {\it Young diagram}
which is a collection of nodes (or boxes) in left-justified rows
with $\lambda_k$ nodes in the $k^{\rm th}$ row.

A partition $\lambda=(\lambda_k)_{k\geq 1}$ is called an {\it
$(m,n)$-hook partition} if $\lambda_{m+1}\leq n$. We denote by
$\mathcal{P}_{m|n}$ the set of all $(m,n)$-hook partitions. A
tableau $T$ obtained by filling a Young diagram $\lambda$ with the
entries in $\B_{m|n}$ is called {\it $(m,n)$-hook semistandard} if
\begin{itemize}
\item[(1)] the entries in each row (resp. column) are weakly increasing from left to right (resp. from top to bottom),
\item[(2)] the entries in $\B_{m|n}^+$ (resp. $\B_{m|n}^-$) are strictly increasing in each column (resp. row)
\end{itemize}
(see \cite{BR}). We say that $\lambda$ is the {\it shape of $T$}. It
is easy to see that a partition $\lambda$ can be made into an
$(m,n)$-hook semistandard tableau if and only if
$\lambda\in\mathcal{P}_{m|n}$.

For $\lambda\in\mathcal{P}_{m|n}$, let $\B_{m|n}(\lambda)$ be the
set of all $(m,n)$-hook semistandard tableaux of shape $\lambda$. We
may view $\B_{m|n}(\lambda)$  as a subset of $\mathcal{W}_{m|n}$ by
column reading (or far eastern reading). That is, we read the
entries of a tableau column by column from right to left, and in
each column we read the entries from top to bottom.

For $T\in \B_{m|n}(\lambda)$, the weight of $T$ is given by
$\text{wt}\,T = \sum_{b\in \B_{m|n}} \mu_b \epsilon_b \in P_{m|n}$,
where $\mu_b$ is number of occurrences of $b$ in $T$. Indeed,
$\B_{m|n}(\lambda)$ together with $0$ is stable under $e_i,f_i$
($i\in I_{m|n}$), and $\B_{m|n}(\lambda)$ is a subcrystal of
$\mathcal{W}_{m|n}$.
\begin{thm}[\cite{BKK}] For $\lambda\in\mathcal{P}_{m|n}$,
$\B_{m|n}(\lambda)$ is a crystal graph for $\frak{gl}_{m|n}$.
Moreover, $\B_{m|n}(\lambda)$ is a connected $I_{m|n}$-colored
oriented graph with a unique highest weight element
$H_{m|n}^{\lambda}$.\qed
\end{thm}

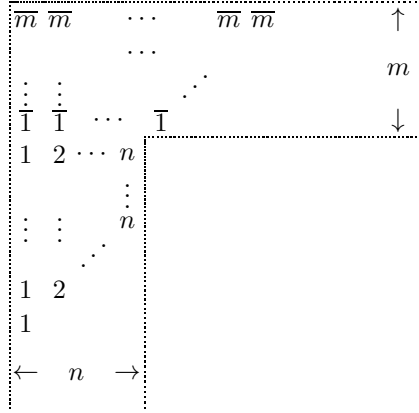
\begin{figure}\label{hwtab}
\setlength{\unitlength}{0.45cm}
\begin{picture}(20,13)(-2,0)
\put(9.7,8.7){\makebox(1,1){.}} \put(10,9){\makebox(1,1){.}}
\put(10.3,9.3){\makebox(1,1){.}}

\qbezier[64](9,8)(13,8)(17,8) \qbezier[96](5,12)(9,12)(17,12)
\put(16,11){\makebox(1,1){$\uparrow$}}
\put(16,9.5){\makebox(1,1){$m$}}
\put(16,8){\makebox(1,1){$\downarrow$}}

\qbezier[96](5,12)(5,6)(5,0) \qbezier[64](9,8)(9,5)(9,0)
\put(5,0.5){\makebox(1,1){$\leftarrow$}}
\put(6.5,0.5){\makebox(1,1){$n$}}
\put(8,0.5){\makebox(1,1){$\rightarrow$}}


\put(7,4){\makebox(1,1){.}}
\put(6.7,3.7){\makebox(1,1){.}} \put(7.3,4.3){\makebox(1,1){.}}

\put(5,11){\makebox(1,1){$\ov{m}$}} \put(6,11){\makebox(1,1){$\ov{m}$}}
\put(11,11){\makebox(1,1){$\ov{m}$}} \put(12,11){\makebox(1,1){$\ov{m}$}}


\put(5,9){\makebox(1,1){$\vdots$}}
\put(6,9){\makebox(1,1){$\vdots$}}
\put(8.5,11){\makebox(1,1){$\cdots$}}
\put(8.5,10){\makebox(1,1){$\cdots$}}

\put(5,8){\makebox(1,1){$\ov{1}$}} \put(6,8){\makebox(1,1){$\ov{1}$}}
\put(7.5,8){\makebox(1,1){$\cdots$}} \put(9,8){\makebox(1,1){$\ov{1}$}}

\put(5,7){\makebox(1,1){${1}$}}
\put(6,7){\makebox(1,1){${2}$}}

\put(5,5){\makebox(1,1){$\vdots$}}
\put(6,5){\makebox(1,1){$\vdots$}}
\put(5,3){\makebox(1,1){${1}$}}
\put(6,3){\makebox(1,1){${2}$}}
\put(5,2){\makebox(1,1){${1}$}}

\put(7,7){\makebox(1,1){$\cdots$}}
\put(8,7){\makebox(1,1){${n}$}}
\put(8,6){\makebox(1,1){$\vdots$}}
\put(8,5){\makebox(1,1){${n}$}}
\end{picture}\ \
\caption{A highest weight tableau $H_{m|n}^{\lambda}$}
\end{figure}

\begin{rem}{\rm
Note that ${\rm wt}(H_{m|n}^{\lambda})\geq {\rm wt}(T)$ for all
$T\in \B_{m|n}(\lambda)$ (\cite{BKK}) (see Figure 1), and hence
$\re_iH^{\lambda}_{m|n}=0$ for all $i\in I_{m|n}$. But, unlike the
crystal graphs associated to integrable highest weight
representations of a symmetrizable Kac-Moody algebra, there might
exist $T\in\B_{m|n}(\lambda)$ such that $T\neq H_{m|n}^{\lambda}$
and $\re_iT=0$ for all $i\in I_{m|n}$. Such a tableau $T$ was called
a {\it fake highest weight vector}, and $H_{m|n}^{\lambda}$ a {\it
genuine highest weight vector} in \cite{BKK}.}
\end{rem}

To characterize a  connected component in $\mathcal{W}_{m|n}$, we
need the algorithm of {\it Schensted's column bumping} for
$(m,n)$-hook semistandard tableaux (\cite{BR,Rem}): for
$\lambda\in\mathcal{P}_{m|n}$ and $T\in\B_{m|n}(\lambda)$, we define
$T \leftarrow{b}$ ($b\in \B_{m|n}$) to be the tableau obtained from
$T$ by applying the following procedure;
\begin{itemize}
\item[(1)] If $b\in \B_{m|n}^+$, let $b'$ be the smallest entry in the first
(or the left-most) column which is greater than or equal to $b$. If
$b\in \B_{m|n}^-$, let $b'$ be the smallest entry in the first
column which is greater than $b$. If there are more than one $b'$,
choose the one in the highest position.

\item[(2)] Replace {$b'$} by {$b$} ({$b'$} is {\it
bumped out} of the first column). If there is no such $b'$, put
{$b$} at the bottom of the first column and stop the procedure.

\item[(3)] Repeat (1) and (2) on the next column with {$b'$}.
\end{itemize}\vskip 3mm

Note that $(T\leftarrow b)\in \B_{m|n}(\mu)$ for some
$\mu\in\mathcal{P}_{m|n}$, where $\mu$ is given by adding a node at
$\lambda$. Now, for a given word $w=w_1\cdots w_r\in\mathcal{W}$, we
define
\begin{equation}\label{column insertion}
P(w)=(\cdots(({w_1}\leftarrow{w_2}\,)\leftarrow{w_3}\,)
\cdots)\leftarrow{w_r}\ .
\end{equation}

\begin{lem}[\cite{BKK}]
For $w\in\mathcal{W}_{m|n}$, we have $w\simeq P(w)$. Hence, any
connected component in $\mathcal{W}_{m|n}$ is isomorphic to
$\B_{m|n}(\lambda)$ for some $\lambda\in \mathcal{P}_{m|n}$.\qed
\end{lem}

We will also use the following lemma in the next section.

\begin{lem}[\cite{KK}]\label{equivlence of tableaux}
Let $T$ and $T'$ be two $(m,n)$-hook semistandard tableaux. If
$T\simeq T'$, then $T=T'$.\qed
\end{lem}

Let $x=\{\,x_b\,|\, b\in \B_{m|n}\,\}$  be the set of variables
indexed by $\B_{m|n}$. For $\mu=\sum_{b\in
\B_{m|n}}\mu_b\epsilon_b\in P_{m|n}$, we set $x^{\mu}=\prod_{b\in
\B_{m|n}}x_b^{\mu_b}$. For $\lambda\in\mathcal{P}_{m|n}$, we define
a {\it hook Schur function} corresponding to $\lambda$ by
\begin{equation*}
hs_{\lambda}(x)=\sum_{T\in\B_{m|n}(\lambda)}x^{{\rm wt}T},
\end{equation*}
which is the character of $\B_{m|n}(\lambda)$ (see \cite{BR,Rem}).

\section{Bicrystal graphs  and Cauchy decomposition}
We consider a bicrystal graph structure on the set of certain
matrices of non-negative integers, which parameterizes the monomial
basis of $S(\mathbb{C}^{m|n}\otimes \mathbb{C}^{u|v})$, and then
derive an explicit decomposition by finding all the highest weight
elements.
\subsection{Crystal graphs of biwords}
Suppose that  $m,n,u,v$ are non-negative integers such that $m+n,
u+v>0$.

Let
\begin{equation}
\begin{split}
\Omega_{m|n,u|v}=\{\,& (\bi,\bj)\in
\mathcal{W}_{m|n}\times \mathcal{W}_{u|v}\,| \\
& \text{(1) $\bi=i_1,\cdots,i_r$ and $\bj=j_1,\cdots,j_r$  for some $r\geq 0$}, \\
& \text{(2) $(i_1,j_1)\leq \cdots \leq (i_r,j_r)$}, \\
& \text{(3) $|i_k|\neq |j_k|$ implies $(i_k,j_k)\neq
(i_{k\pm1},j_{k\pm 1})$}\},
\end{split}
\end{equation}
where for $(i,j)$ and $(k,l)\in \B_{m|n}\times \B_{u|v}$, the {\it
super lexicographic ordering} is given by
\begin{equation}\label{partial order}
(i,j)< (k,l) \ \ \ \ \Leftrightarrow \ \ \ \
\begin{cases}
(j<l) & \text{or}, \\
(j=l\in \B_{u|v}^+,\ \text{and} \ i>k) & \text{or}, \\
(j=l\in \B_{u|v}^-,\ \text{and} \ i<k) &.
\end{cases}
\end{equation}
Next, let $\Omega^*_{m|n,u|v}$ be the set of pairs $(\bk,\bl)\in
\mathcal{W}_{m|n}\times \mathcal{W}_{u|v}$ such that $(\bl,\bk)\in
\Omega_{u|v,m|n}$. For simplicity, we write
$\Omega=\Omega_{m|n,u|v}$ and $\Omega^*=\Omega^*_{m|n,u|v}$.

Now, for $i\in I_{m|n}$ and $(\bi,\bj)\in \Omega$, we define
\begin{equation*}
e_i (\bi,\bj)=(\re_i\bi,\bj), \ \ \ \ \ f_i
(\bi,\bj)=(\rf_i\bi,\bj),
\end{equation*}
where we assume that $x_i (\bi,\bj)=0$ if $x_i\bi=0$ ($x=e,f$). We
set ${\rm wt}(\bi,\bj)={\rm wt}(\bi)$,
$\varepsilon_i(\bi,\bj)=\varepsilon_i(\bi)$ and
$\varphi_i(\bi,\bj)=\varphi_i(\bi)$ ($i\in I_{m|n}$).

 Similarly,
for $i\in I_{u|v}$ and $(\bk,\bl)\in \Omega^*$, we define
\begin{equation*}
e_j^* (\bk,\bl)=(\bk,\re_j\bl), \ \ \ \ \ f_j^*
(\bk,\bl)=(\bk,\rf_j\bl),
\end{equation*}
and set ${\rm wt}^*(\bk,\bl)={\rm wt}(\bl)$,
$\varepsilon^*_j(\bk,\bl)=\varepsilon_j(\bl)$ and
$\varphi^*_j(\bk,\bl)=\varphi_j(\bl)$ ($j\in I_{u|v}$).

\begin{lem}\label{crystal} Under the above hypothesis,
\begin{itemize}
\item[(1)] the set $\Omega$ together with ${\rm wt}, e_i,f_i,\varepsilon_i,\varphi_i$ $(i\in I_{m|n})$
is a crystal graph for $\frak{gl}_{m|n}$,

\item[(2)] the set $\Omega^*$ together with ${\rm wt}^*, e^*_j,f^*_j,\varepsilon^*_j,\varphi^*_j$ $(j\in I_{u|v})$
is a crystal graph for $\frak{gl}_{u|v}$.
\end{itemize}
\end{lem}
\pf We will prove only (1) since the proof of (2) is the same.
Suppose that $(\bi,\bj)\in\Omega\setminus \{(\emptyset,\emptyset)\}$
is given where $\bi=i_1\cdots i_r$ and $\bj=j_1\cdots j_r$ for some
$r\geq 1$. We write
$$\bi=\bi_{\ov{u}}\bi_{\ov{u-1}}\cdots\bi_{v-1} \bi_v,$$ where
$\bi_b=i_{t_1}\cdots i_{t_b}$  ($b\in \B_{u|v}$) is a subword of
$\bi$ such that $j_{t_1}=\cdots=j_{t_b}=b$.

Then, $\bi_b$ ($b\in\B_{u|v}$) is $\frak{gl}_{m|n}$-equivalent to an
$(m,n)$-hook semistandard tableau $T_b$ of a single row or a single
column as follows;
\begin{equation*}
T_b=
\begin{cases}
\begin{array}{ccc}
     i_{t_1}\cdots i_{t_b}
\end{array} \in \B_{m|n}((t_b)), & \text{if $b\in\B_{u|v}^+$}, \\
\ \ \ \begin{array}{c}
  i_1 \\
  \vdots \\
  i_{t_b} \\
\end{array} \ \ \ \ \in \B_{m|n}((1^{t_b})),
& \text{if $b\in\B_{u|v}^-$}.
\end{cases}
\end{equation*}
So, $\bi$ is $\frak{gl}_{m|n}$-equivalent to
\begin{equation*}
T_{\ov{u}}\otimes T_{\ov{u-1}}\otimes \cdots\otimes T_{v-1}\otimes
T_v,
\end{equation*}
where $T_b=\emptyset$ if $\bi_b=\emptyset$. Therefore, if
$x_i(\bi,\bj)\neq 0$ ($x=e,f$) for $i\in I_{m|n}$, then
$(x_i\bi,\bj)\in\Omega$. It follows that $\Omega$ is a crystal graph
for $\frak{gl}_{m|n}$. \qed

\subsection{Bicrystal graphs}
Consider the following set of matrices of non-negative integers;
\begin{equation}
\begin{split}
\M_{m|n,u|v}= \{\,A=&(a_{bb'})_{b\in \B_{m|n}, b'\in \B_{u|v}}\,|\, \\
& \text{(1) $a_{bb'}\in\mathbb{Z}_{\geq 0}$,\ \  (2) $a_{bb'}\leq 1$
if $|b|\neq|b'|$}\, \}.
\end{split}
\end{equation}
For simplicity, we write $\M=\M_{m|n,u|v}$.

For $(\bi,\bj)\in \Omega$, define $A(\bi,\bj)=(a_{bb'})$ to be the
matrix in $\M$, where $a_{bb'}$ is the number of $k$'s such that
$(i_k,j_k)=(b,b')$ for $b\in \B_{m|n}$ and $b'\in \B_{u|v}$. Then,
it follows that the map $(\bi,\bj)\mapsto A(\bi,\bj)$ gives a
bijection between $\Omega$ and $\M$, where the pair of empty words
$(\emptyset,\emptyset)$ corresponds to zero matrix. Similarly, we
have a bijection $(\bk,\bl)\mapsto A(\bk,\bl)$ from $\Omega^*$ to
$\M$.

With these identifications, $\M$ becomes a crystal graph for both
$\frak{gl}_{m|n}$ and $\frak{gl}_{u|v}$ by Lemma \ref{crystal}.

\begin{ex}\label{example1}{\rm
Suppose that $m|n=u|v=2|2$ and
\begin{equation*}
A=\left(\begin{tabular}{cc|cc}
  1 & 0 & 1 & 1 \\
  0 & 2 & 0 & 0 \\ \hline
  1 & 0 & 0 & 1 \\
  0 & 0 & 2 & 0 \\
\end{tabular}\right )\in\M.
\end{equation*}
Then $A=A(\bi,\bj)=A(\bk,\bl)$ for $(\bi,\bj)\in\Omega$ and
$(\bk,\bl)\in\Omega^*$, where
\begin{equation*}
\begin{split}
&\bi=1 \ \ov{2} \ \ov{1} \ \ov{1} \ \ov{2} \ 2 \ 2 \ \ov{2} \ 1, \ \ \ \
\bk=\ov{2} \ \ov{2} \ \ov{2} \ \ov{1} \ \ov{1} \ 1 \ 1 \ 2 \ 2,  \\
&\bj=\ov{2} \ \ov{2} \ \ov{1} \ \ov{1} \ 1 \ 1 \ 1 \ 2 \ 2, \ \ \ \
\bl\,\, =2 \ 1 \ \ov{2} \ \ov{1} \ \ov{1} \ \ov{2} \ 2 \ 1 \ 1.
\end{split}
\end{equation*}
}
\end{ex}

\begin{df}{\rm
Let $B$ be a crystal graph for both $\frak{gl}_{m|n}$ and
$\frak{gl}_{u|v}$. We denote by $e_j^*$ and $f_j^*$ ($j\in I_{u|v}$)
the Kashiwara operators for $\frak{gl}_{u|v}$. We call $B$ a {\it
crystal graph for $\frak{gl}_{m|n}\oplus\frak{gl}_{u|v}$, or
$(\frak{gl}_{m|n},\frak{gl}_{u|v})$-bicrystal} if $e_i,f_i$ commute
with $e_j^*,f_j^*$ ($i\in I_{m|n},j\in I_{u|v}$), where we
understand the Kashiwara operators as the associated maps from
$B\cup\{0\}$ to itself (that is, $x_i0=x_j^*0=0$, for $x=e,f$). }
\end{df}
 For example,
$B=\B_{m|n}(\lambda)\times \B_{u|v}(\mu)$ ($\lambda\in
\mathcal{P}_{m|n}$, $\mu\in \mathcal{P}_{u|v}$) is a
$(\frak{gl}_{m|n},\frak{gl}_{u|v})$-bicrystal where
$x_i(b_1,b_2)=(x_ib_1,b_2)$, $x^*_j(b_1,b_2)=(b_1,x^*_jb_2)$ for
$(b_1,b_2)\in B$ and $x=e,f$.

\begin{lem}\label{bicrystal}\mbox{}
$\M$ is a $(\frak{gl}_{m|n},\frak{gl}_{u|v})$-bicrystal.
\end{lem}
\pf The proof is a straightforward verification. So, let us prove
the following case;
\begin{equation}\label{commuting}
e_t^*e_s A =e_se_t^*A
\end{equation}
for $A\in\M$, $1\leq s\leq n-1$ and $1\leq t\leq v-1$. The other
cases can be checked in a similar manner.

To show this, we may assume that $A=(a_{bb'})$ such that $a_{bb'}=0$
unless $b\in\{s,s+1\}$ or $b'\in\{t,t+1\}$. Put
$A'=(a_{bb'})_{b=s,s+1}$ and $A''=(a_{bb'})_{b'=t,t+1}$, which are
the submatrices of $A$. If neither $A'$ nor $A''$ contains both an
$s$-good $-$ sign and a $t$-good $-$ sign of $A$, then
\eqref{commuting} follows from the fact that the $s$-signature of
$e_t^*A$ is the same as that of $A$, and the $t$-signature of $A$ is
the same as that of $e_sA$. So we may assume that $A$ is either $A'$
or $A''$, and prove the case only when $A=A'$, that is,
\begin{equation*}
A=\left(%
\begin{array}{cccccc}
  x^+_{\ov{u}} & \cdots & x^+_{\ov{1}} & x^+_{1} & \cdots & x^+_{v} \\
  x^-_{\ov{u}} & \cdots & x^-_{\ov{1}} & x^-_{1} & \cdots & x^-_{v} \\
\end{array}%
\right),
\end{equation*}
where $x^+_b=a_{sb}$ and $x^-_b=a_{s+1 b}$ ($b\in\B_{u|v}$). Note
that $A=A(\bi,\bj)=A(\bk,\bl)$ for unique $(\bi,\bj)\in\Omega$ and
$(\bk,\bl)\in\Omega^*$. Then $\bi$ and $\epsilon^{(s)}(\bi)$ are
of the form:
\begin{equation*}
\begin{split}
\bi&=(s+1)^{x^-_{\ov{u}}}\, s^{x^+_{\ov{u}}}\,\cdots\,
(s+1)^{x^-_{\ov{1}}}\, s^{x^+_{\ov{1}}}\,
s^{x^+_{1}}\,(s+1)^{x^-_{1}}\,\cdots\, s^{x^+_{v}}\,(s+1)^{x^-_{v}},\\
\epsilon^{(s)}(\bi)&=(-^{x^-_{\ov{u}}},+^{x^+_{\ov{u}}},\cdots,
-^{x^-_{\ov{1}}},+^{x^+_{\ov{1}}}, +^{x^+_{1}},-^{x^-_{1}},\cdots,
+^{x^+_{v}},-^{x^-_{v}}),
\end{split}
\end{equation*}
where the multiplicities of letters and signs are given as
exponents. On the other hand, the $t$-signature of $A$ or $\bl$ is
completely determined by its subword $\bl'$:
\begin{equation*}
\begin{split}
\bl'&=t^{x^+_t}\, (t+1)^{x^+_{t+1}}\, t^{x^-_t}\,
(t+1)^{x^-_{t+1}}, \\
\epsilon^{(t)}(\bl')&=(+^{x^+_t}, -^{x^+_{t+1}}, +^{x^-_t},
-^{x^-_{t+1}}).
\end{split}
\end{equation*}
For convenience, we set
\begin{equation*}
\tilde{A}=\left(%
\begin{array}{cc}
  a & b \\
  c & d \\
\end{array}%
\right)=
\left(%
\begin{array}{cc}
  x^+_{t} & x^+_{t+1} \\
  x^-_{t} & x^-_{t+1} \\
\end{array}%
\right).
\end{equation*}

\textsc{Case 1}. Suppose that $e_t^*A=0$ (equivalently,
$e_t^*\tilde{A}=0$). This implies that $d=0$ and $b\leq c$ (note
that in this case, we cancel out $(-,+)$ pairs to obtain a
$t$-signature). If an $s$-good $-$ sign occurs in $\tilde{A}$, then
we have $b<c$ and
\begin{equation*}
e_sA=\left(%
\begin{array}{cccc}
  \cdots & a+1 & b & \cdots  \\
  \cdots & c-1 & 0 & \cdots  \\
\end{array}%
\right).
\end{equation*}
Since $b\leq c-1$, we get $e_t^*e_sA=0$, which implies
\eqref{commuting}. If an $s$-good $-$ sign of $A$ does not occur in
$\tilde{A}$, then $\tilde{A}$ does not change when we apply $e_s$ to
$A$ and $e_t^*e_sA=0=e_se_t^*A$, which also implies
\eqref{commuting}.\vskip 3mm

\textsc{Case 2}. Suppose that $e_t^*A\neq 0$ and $b$ is changed by
$e_t^*$. This implies that $b>0$ and $b>c$. If
$e_t^*A=A(\bi',\bj')$ for $(\bi',\bj')\in\Omega$, then
\begin{equation*}
\begin{split}
\epsilon^{(s)}(\bi)&=(\cdots,+^{a},-^c,+^{b},-^d,\cdots), \\
\epsilon^{(s)}(\bi')&=(\cdots,+^{a+1},-^c,+^{b-1},-^d,\cdots).
\end{split}
\end{equation*}
Note that the subsequences $(+^{a+1},-^c,+^{b-1},-^d)$ and
$(+^{a},-^c,+^{b},-^d)$ reduce to the same sequence
$(+^{a+b-c},-^d)$ ($a+b-c>0$), and there is no $s$-good $-$ sign in
$c$. If no $s$-good $-$ sign of $A$ occurs in $\tilde{A}$ (or, in
$d$), then it is easy to see that \eqref{commuting} holds. If there
is an $s$-good $-$ sign in $\tilde{A}$, then
\begin{equation*}
\begin{split}
&e_se_t^*\left(%
\begin{array}{cc}
    a & b   \\
    c & d   \\
\end{array}%
\right)=e_s
\left(%
\begin{array}{cc}
   a+1 & b-1   \\
   c & d   \\
\end{array}%
\right)=
\left(%
\begin{array}{cccc}
   a+1 & b  \\
   c & d-1   \\
\end{array}%
\right),\\
&e_t^*e_s\left(%
\begin{array}{cc}
    a & b   \\
    c & d   \\
\end{array}%
\right)=e_t^*
\left(%
\begin{array}{cc}
   a & b+1   \\
   c & d-1   \\
\end{array}%
\right)=
\left(%
\begin{array}{cccc}
   a+1 & b  \\
   c & d-1   \\
\end{array}%
\right).
\end{split}
\end{equation*}
Hence, we have \eqref{commuting}.\vskip 3mm

\textsc{Case 3}. Suppose that $e_t^*A\neq 0$ and $d$ is changed by
$e_t^*$. This implies that $d>0$ and $b\leq c$. If
$e_t^*A=A(\bi',\bj')$ for $(\bi',\bj')\in\Omega$, then
\begin{equation*}\label{t-signature}
\begin{split}
\epsilon^{(s)}(\bi)&=(\cdots,+^{a},-^c,+^{b},-^d,\cdots), \\
\epsilon^{(s)}(\bi')&=(\cdots,+^{a},-^{c+1},+^b,-^{d-1},\cdots).
\end{split}
\end{equation*}
Note that the subsequences $(+^{a},-^{c+1},+^b,-^{d-1})$ and
$(+^{a},-^c,+^{b},-^d)$ reduce to the same sequence
$(+^{a},-^{c-b+d})$ ($c-b+d>0$). If no $s$-good $-$ sign of $A$
occurs in $\tilde{A}$ (or, in $d$), then it is easy to see that
$e_t^*e_sA=e_se_t^*A$. Assume that an $s$-good $-$ sign of $A$
occurs in $\tilde{A}$. If $b=c$, then we have
\begin{equation*}
\begin{split}
&e_se_t^*\left(%
\begin{array}{cc}
    a & b   \\
    b & d   \\
\end{array}%
\right)=e_s
\left(%
\begin{array}{cc}
   a & b   \\
   b+1 & d-1   \\
\end{array}%
\right)=
\left(%
\begin{array}{cccc}
   a+1 & b  \\
   b & d-1   \\
\end{array}%
\right),\\
&e_t^*e_s\left(%
\begin{array}{cc}
    a & b   \\
    b & d   \\
\end{array}%
\right)=e_t^*
\left(%
\begin{array}{cc}
   a & b+1   \\
   b & d-1   \\
\end{array}%
\right)=
\left(%
\begin{array}{cccc}
   a+1 & b  \\
   b & d-1   \\
\end{array}%
\right).
\end{split}
\end{equation*}
If $b<c$, then we have
\begin{equation*}
\begin{split}
&e_se_t^*\left(%
\begin{array}{cc}
    a & b   \\
    c & d   \\
\end{array}%
\right)=e_s
\left(%
\begin{array}{cc}
   a & b   \\
   c+1 & d-1   \\
\end{array}%
\right)=
\left(%
\begin{array}{cccc}
   a+1 & b  \\
   c & d-1   \\
\end{array}%
\right),\\
&e_t^*e_s\left(%
\begin{array}{cc}
    a & b   \\
    c & d   \\
\end{array}%
\right)=e_t^*
\left(%
\begin{array}{cc}
   a+1 & b   \\
   c-1 & d   \\
\end{array}%
\right)=
\left(%
\begin{array}{cccc}
   a+1 & b  \\
   c & d-1   \\
\end{array}%
\right).
\end{split}
\end{equation*}
Hence, we have \eqref{commuting}. \qed \vskip 3mm

\begin{lem}\label{connected component} \mbox{}
\begin{itemize}
\item[(1)] Let $C$ be a connected component in $\M$ as a
$\frak{gl}_{m|n}$-crystal. If $x_j^*C\neq \{0\}$ for some $j\in
I_{u|v}$ and $x=e,f$, then $x_j^* : C \rightarrow x_j^*C$ is an
isomorphism of $\frak{gl}_{m|n}$-crystals.

\item[(2)] Let $C^*$ be a connected component in $\M$ as a
$\frak{gl}_{u|v}$-crystal. If $x_iC^*\neq \{0\}$ for some $i\in
I_{m|n}$ and $x=e,f$, then $x_i : C^* \rightarrow x_iC^*$ is an
isomorphism of $\frak{gl}_{u|v}$-crystals.
\end{itemize}
\end{lem}
\pf (1)  Choose $b\in C$ such that $x_j^*b\neq 0$. For any $b'\in
C$, we have
\begin{equation*}
b'=x_{i_1}\cdots x_{i_t}b,
\end{equation*}
for some Kashiwara operators $x_{i_k}$  ($i_k\in I_{m|n}, 1\leq
k\leq t$). We first claim that $x_j^*b'\neq 0$.

We will use induction on $t$. Suppose that $t=1$. Let $y_{i_1}$ be
the Kashiwara operator given by $y_{i_1}=f_{i_1}$ (resp.
$y_{i_1}=e_{i_1}$) if $x_{i_1}=e_{i_1}$ (resp. $x_{i_1}=f_{i_1}$).
Then $b=y_{i_1}b'$, and $$0\neq
x_j^*b=x_j^*y_{i_1}b'=y_{i_1}x_j^*b'$$ by Lemma \ref{bicrystal},
which implies that $x_j^*b'\neq 0$. Suppose that $t>1$. Since
$x_j^*x_{i_t}b\neq 0$ and $b'=x_{i_1}\cdots x_{i_{t-1}}(x_{i_t}b)$,
it follows that $x_j^*b'\neq 0$ by induction hypothesis. This
completes the induction.

Hence, the composite of the following two maps is the identity map
on $C$;
\begin{equation*}
C \stackrel{x_j^*}{\longrightarrow} x_j^*C
\stackrel{y_j^*}{\longrightarrow} C,
\end{equation*}
where $y_j^*=e_j^*$ (resp. $f_j^*$) if $x_j^*=f_j^*$ (resp.
$e_j^*$). This implies that $x_j^*$ is a bijection which commutes
with $e_i,f_i$ ($i\in I_{m|n}$) by Lemma \ref{bicrystal}, and that
$x_j^*C$ is isomorphic to $C$ as a $\frak{gl}_{m|n}$-crystal. The
proof of (2) is similar. \qed \vskip 3mm

Given $A\in \M$, suppose that $A=A(\bi,\bj)=A(\bk,\bl)$ for unique
$(\bi,\bj)\in \Omega$ and $(\bk,\bl)\in \Omega^*$. We define
\begin{equation}\label{pi}
\pi(A)=(P_1(A),P_2(A))=(P(\bi),P(\bl)),
\end{equation}
(see \eqref{column insertion}). Then
$\pi(A)\in\B_{m|n}(\lambda)\times \B_{u|v}(\mu)$ for some
$(\lambda,\mu)\in\mathcal{P}_{m|n}\times\mathcal{P}_{u|v}$. Note
that $A\simeq_{\frak{gl}_{m|n}}P_1(A)$, and
$A\simeq_{\frak{gl}_{u|v}}P_2(A)$.

\begin{lem}\label{P-tableaux}
Suppose that $A\in\M$ is given.
\begin{itemize}
\item[(1)] if $x_j^*A\neq 0$ for some $j\in I_{u|v}$ and $x=e,f$,
then $P_1(x_j^*A)=P_1(A)$.

\item[(2)] if $x_iA\neq 0$ for some $i\in I_{m|n}$ and $x=e,f$,
then $P_2(x_iA)=P_2(A)$.
\end{itemize}
\end{lem}
\pf It suffices to prove (1). By Lemma \ref{connected component},
$A$ is $\frak{gl}_{m|n}$-equivalent to $x_j^*A$, and hence the
tableaux $P_1(A)$ and $P_1(x_j^*A)$ are
$\frak{gl}_{m|n}$-equivalent. By Lemma \ref{equivlence of tableaux},
we have $P_1(x_j^*A)=P_1(A)$. \qed \vskip 3mm

Let $B$ be a  $(\frak{gl}_{m|n},\frak{gl}_{u|v})$-bicrystal. If we
set ${\mathcal I}=I_{m|n}\sqcup {I^*_{u|v}}$ where
${I^*_{u|v}}=\{j^*|j\in I_{u|v}\}$, then $B$ is an ${\mathcal
I}$-colored oriented graph with respect to $x_i,x_j^*$ for $i\in
I_{m|n},j\in I_{u|v}$, and $x=e,f$. Now, we can characterize a
connected component in $\M$ as a
$(\frak{gl}_{m|n},\frak{gl}_{u|v})$-bicrystal.
\begin{prop}\label{g-crystal}
For each connected component $\mathbf{C}$ in $\M$, $\pi$ gives the
following isomorphism of
$(\frak{gl}_{m|n},\frak{gl}_{u|v})$-bicrystals;
$$\pi: \mathbf{C} \longrightarrow \B_{m|n}(\lambda)\times
\B_{u|v}(\mu),$$ for some $\lambda\in \mathcal{P}_{m|n}$ and $\mu\in
\mathcal{P}_{u|v}$.
\end{prop}
\pf By Lemma \ref{P-tableaux}, we have
\begin{equation*}
\pi(x_iA)=(x_iP_1(A),P_2(A)),
\end{equation*}
for $A\in \M$, $i\in I_{m|n}$ and $x=e,f$ (we assume that $\pi(0)=0$
and $(0,P_2(A))=0$). Similarly, we have
$\pi(x_j^*A)=(P_1(A),x_jP_2(A))$ for $j\in I_{u|v}$.

Let $\mathbf{C}$ be a connected component in $\M$ as an
$\mathcal{I}$-colored oriented graph. Choose an arbitrary $A\in
\mathbf{C}$. Suppose that $\pi(A)\in\B_{m|n}(\lambda)\times
\B_{u|v}(\mu)$ for some $\lambda\in\mathcal{P}_{m|n}$ and
$\mu\in\mathcal{P}_{u|v}$. Then we have
\begin{equation*}
\pi: \mathbf{C} \longrightarrow \B_{m|n}(\lambda)\times
\B_{u|v}(\mu),
\end{equation*}
where $\pi$ commutes with $e_i,f_i$ and $e_j^*,f_j^*$ $(i\in
I_{m|n}, j\in I_{u|v})$. It is clear that $\pi$ is onto.

Now suppose that $\pi(A)=\pi(A')$ for some $A,A'\in \mathbf{C}$.
Since $\mathbf{C}$ is connected,
$$x_{j_1}^*\cdots x_{j_q}^*x_{i_1}\cdots x_{i_p}A=A',$$
for some  Kashiwara operators $x_{i_k}$ ($i_k\in I_{m|n}$, $1\leq
k\leq p$), $x^*_{j_l}$ ($j_l\in I_{u|v}$, $1\leq l\leq q$). Put
$A''=x_{i_1}\cdots x_{i_p}A$. Then $A$ and $A''$ belong to the same
connected component as a $\frak{gl}_{m|n}$-crystal, say $C_1$. On
the other hand, by Lemma \ref{P-tableaux}, we have
$$P_1(A'')=P_1(A')=P_1(A).$$
Since the map $P_1 : C_1 \rightarrow \B_{m|n}(\lambda)$ is an
isomorphism of $\frak{gl}_{m|n}$-crystals, it follows that $A''=A$.

Next, $A''$ and $A'$ belong to the same connected component as a
$\frak{gl}_{u|v}$-crystal, say $C_2$. By Lemma \ref{P-tableaux}, we
have
$$P_2(A'')=P_2(A)=P_2(A').$$
Since the map $P_2 : C_2 \rightarrow \B_{u|v}(\mu)$ is an
isomorphism, we have $A''=A'$, and hence $A=A'$. So, $\pi$ is
one-to-one.

Therefore, $\pi$ is an isomorphism of
$(\frak{gl}_{m|n},\frak{gl}_{u|v})$-bicrystals. \qed

\begin{ex}{\rm
Let $A$ be as in Example \ref{example1}. Then
\begin{equation*}
\pi(A)=\left(\begin{array}{ccccc}
\ov{2} & \ov{2} & \ov{2} & \ov{1} &  1 \\
\ov{1} & 2 &    &  &  \\
1 &  &  &  &   \\
2 &  &  &  &
\end{array},
\begin{array}{ccccc}
\ov{2} & \ov{2} & \ov{1} & 1 &  2 \\
\ov{1} & 2 &  &  &    \\
1 &  &  &  &    \\
2 &  &  &  &
\end{array}\right ),
\end{equation*}
and the connected component of $A$ is isomorphic to
$\B_{2|2}(5,2,1,1)\times \B_{2|2}(5,2,1,1)$.
 }
\end{ex}
\vskip 3mm

\subsection{Decomposition of $\M$}
Now, we will describe an explicit decomposition of $\M$. Set
\begin{equation}
\M_{\rm
h.w.}=\{\,A\in\M\,|\,\pi(A)=(H_{m|n}^{\lambda},H^{\mu}_{u|v})\
\text{for some
$(\lambda,\mu)\in\mathcal{P}_{m|n}\times\mathcal{P}_{u|v}$}\,\},
\end{equation}
which is the set of all the highest weight elements in $\M$. We have
seen in the proof of Proposition \ref{g-crystal} that $\pi$ induces
an isomorphism of $(\frak{gl}_{m|n},\frak{gl}_{u|v})$-bicrystals
between the connected components of $A\in \M$ and $\pi(A)$. Hence,
\begin{equation*}
\M=\bigoplus_{A\in\M_{\rm h.w.}}\mathbf{C}(A),
\end{equation*}
where $\mathbf{C}(A)$ is the connected component of $A$ in $\M$,
which is isomorphic to $\B_{m|n}(\lambda)\times \B_{u|v}(\mu)$ for
some $\lambda\in\mathcal{P}_{m|n}$ and $\mu\in\mathcal{P}_{u|v}$.

Suppose that
$\lambda=(\lambda_1,\cdots,\lambda_r)\in\mathcal{P}_{m|n}\cap\mathcal{P}_{u|v}$
is given. Let $\nu=(\nu_1,\nu_2,\cdots,\nu_{\ell})$ be the sequence
of non-negative integers ($\ell=\lambda_{m+1}$) determined by
\begin{equation}\label{nu}
(\nu_1+\cdots+\nu_{\ell},\nu_2+\cdots+\nu_{\ell},\cdots,\nu_{\ell})=(\lambda_{m+1},\cdots,\lambda_r)',
\end{equation}
where $(\lambda_{m+1},\cdots,\lambda_r)'$ is the conjugate of the
partition $(\lambda_{m+1},\cdots,\lambda_r)$ (cf.\cite{Mac}).

Assume that $m\geq u$. Let us define $A_{\lambda}=(a_{bb'})\in \M$
by
\begin{itemize}
\item[(1)] for $0\leq k< m$, $0\leq l<u$,
\begin{equation*}
a_{\ov{m-k}\, \ov{u-l}}=
\begin{cases}
\lambda_{k+1}, & \text{if $0\leq k=l< u$}, \\
0, & \text{otherwise}.
\end{cases}
\end{equation*}

\item[(2)] for $0\leq k< m$, $1\leq t\leq v$,
\begin{equation*}
a_{\ov{m-k}\, t}=
\begin{cases}
1, & \text{if $u\leq k< m$ and $\lambda_{k+1}\geq t$}, \\
0, & \text{otherwise}.
\end{cases}
\end{equation*}

\item[(3)] for $1\leq s\leq n$, $0\leq l< u$, \ \ $a_{s\,\ov{u-l}}=0$.

\item[(4)] for $1\leq s\leq n$, $1\leq t\leq v$,
\begin{equation*}
a_{s\, t}=
\begin{cases}
\nu_{s+t-1}, & \text{if $2\leq s+t \leq \ell+1$}, \\
0, & \text{otherwise}.
\end{cases}
\end{equation*}
\end{itemize}

Note that $A_{\lambda}$ is of the following form;
\begin{equation}\label{Alambda}
A_{\lambda}=\left(%
\begin{array}{ccccc|cccccccc}
  \lambda_1 & & & & & & & & & & & & \\
   & \lambda_2 & & & & & & & & & & & \\
   & & \ddots & & & & & & & & & &\\
   &  & &  \lambda_{u-1} & & & & & & & & &\\
   &  & & &  \lambda_{u} & & & & & & & & \\
   &  & & & & 1 & \cdots & \cdots & \cdots & 1 & & & \\
   &  & & & &  \vdots &  & & &  & & &\\
   &  & & & & 1 & \cdots & \cdots  &\ \ \  1 &  & & & \\\hline
   &  & & & &  \nu_{1} & \nu_{2} & \cdots & \nu_{\ell} & & & &\\
   &  & & & & \nu_{2} & \nu_{3} & \cdots &  & & & & \\
   &  & & & & \vdots &  &  &  & & & & \\
   & & & & & \nu_{\ell} &  &  &  & &  & & \\
   &  & & & & &  &  &  &  & & &
\end{array}%
\right).
\end{equation}
If $m\leq u$, then we define $A_{\lambda}$ to be the transpose of
\eqref{Alambda}, where $m$ and $u$ are exchanged. Note that
$A_{\lambda}$ is uniquely determined by $\lambda$. By Schensted's
column bumping algorithm, it is not difficult to see that
\begin{lem} For $\lambda \in \mathcal{P}_{m|n}\cap \mathcal{P}_{u|v}$, we have
\begin{equation*}
\pi(A_{\lambda})=(H_{m|n}^{\lambda},H_{u|v}^{\lambda})\in\B_{m|n}(\lambda)\times\B_{u|v}(\lambda).
\end{equation*}\qed
\end{lem}
\begin{ex}{\rm Suppose that $m|n=4|3$ and $u|v=3|3$. Let $\lambda=(7,5,5,3,3,2,2,1)$. Then
\begin{equation*}
A_{\lambda}=\left(%
\begin{array}{ccc|ccc}
  7 & 0 & 0 & 0 & 0 & 0  \\
  0 & 5 & 0 & 0 & 0 & 0  \\
  0 & 0 & 5 & 0 & 0 & 0  \\
  0 & 0 & 0 & 1 & 1 & 1  \\ \hline
  0 & 0 & 0 & 1 & 2 & 1  \\
  0 & 0 & 0 & 2 & 1 & 0  \\
  0 & 0 & 0 & 1 & 0 & 0
\end{array}%
\right),
\end{equation*}
and
\begin{equation*}
\pi(A_{\lambda})=\left(
\begin{array}{ccccccc}
  \ov{4}  & \ov{4}  & \ov{4}  & \ov{4}  & \ov{4}  & \ov{4} & \ov{4} \\
  \ov{3}  & \ov{3}  & \ov{3}  & \ov{3}  & \ov{3}  &  & \\
  \ov{2}  & \ov{2}  & \ov{2}  & \ov{2}  & \ov{2} &  & \\
  \ov{1}  & \ov{1}  & \ov{1}  &    &   &  & \\
  1  & 2  & 3  &   &   &  & \\
  1  & 2  &   &   &   &  & \\
  1  & 2  &   &   &   &  & \\
  1  &   &   &   &   &  & \\
\end{array},
\begin{array}{ccccccc}
  \ov{3}  & \ov{3}  & \ov{3}  & \ov{3}  & \ov{3}  & \ov{3} & \ov{3} \\
  \ov{2}  & \ov{2}  & \ov{2}  & \ov{2}  & \ov{2}  &  & \\
  \ov{1}  & \ov{1}  & \ov{1}  & \ov{1}  & \ov{1} &  & \\
  1  & 2  & 3  &    &   &  & \\
  1  & 2  & 3  &   &   &  & \\
  1  & 2  &   &   &   &  & \\
  1  & 2  &   &   &   &  & \\
  1  &   &   &   &   &  & \\
\end{array}\right ).
\end{equation*}
}
\end{ex}

\begin{thm}\label{highest weight vector1} We have
\begin{equation*}
\M_{\rm h.w.}=\{\,A_{\lambda}\,|\,\lambda\in
\mathcal{P}_{m|n}\cap\mathcal{P}_{u|v}\,\},
\end{equation*}
and hence the following isomorphism of
$(\frak{gl}_{m|n},\frak{gl}_{u|v})$-bicrystals;
\begin{equation*}
\pi : \M \longrightarrow \bigoplus_{\lambda\in
\mathcal{P}_{m|n}\cap\mathcal{P}_{u|v}} \B_{m|n}(\lambda)\times
\B_{u|v}(\lambda).
\end{equation*}
\end{thm}
\pf  For convenience, we assume that $m\geq u$ and $n\geq v$.
Suppose that $A=(a_{bb'})\in\M_{\rm h.w.}$ is given. We claim that
$A=A_{\lambda}$ for some $\lambda\in
\mathcal{P}_{m|n}\cap\mathcal{P}_{u|v}$ (see \eqref{Alambda}). We
may assume that $A$ is a non-zero matrix since the zero matrix
corresponds to $A_{(0)}$.  \vskip 3mm

\textsc{Step 1.} Suppose that $A=A(\bi,\bj)=A(\bk,\bl)$ for
$(\bi,\bj)\in\Omega$, $(\bk,\bl)\in\Omega^*$, and $\bi=i_1\cdots
i_r$ and $\bj=j_1\cdots j_r$ for some $r\geq 1$. Let us write
$$\bi=\bi_{\ov{u}}\bi_{\ov{u-1}}\cdots\bi_{v-1} \bi_v,$$ where
$\bi_b=i_{t_1}\cdots i_{t_b}$ ($b\in \B_{u|v}$, $t_b\geq 0$) is a
subword of $\bi$ such that $j_{t_1}=\cdots=j_{t_b}=b$. Similarly, we
write $\bl=\bl_{\ov{m}}\cdots \bl_{n}$.

Since $A\in \M_{\rm h.w.}$, $\pi(A)=(P(\bi),P(\bl))$ is a pair of
highest weight tableaux. By Schensted's algorithm, we observe that
the shape of $P(\bi_{\ov{u}})$ is a single  row, and all the letters
in $\bi_{\ov{u}}$ are placed in the first row of $P(\bi)$. Hence, we
should have $\bi_{\ov{u}}=\ov{m}\cdots \ov{m}$ ($a_{\ov{m}\,\ov{u}}$
times), equivalently, $a_{b\ov{u}}=0$ for all $b\in
\B_{m|n}\setminus\{\ov{m}\}$. By the same arguments, we have
$\bl_{\ov{m}}=\ov{u}\cdots \ov{u}$ ($a_{\ov{m}\,\ov{u}}$ times), or
$a_{\ov{m}b}=0$ for all $b\in \B_{u|v}\setminus\{\ov{u}\}$.

Next, consider $\bi_{\ov{u-1}}$. Since $a_{\ov{m}\ov{u-1}}=0$, all
the letters in $\bi_{\ov{u-1}}$ are greater than $\ov{m}$, and
placed in the first two rows of $P(\bi)$. This implies that
$\bi_{\ov{u-1}}=\ov{m-1}\cdots \ov{m-1}$ ($a_{\ov{m-1}\,\ov{u-1}}$
times), and $a_{\ov{m}\,\ov{u}}\geq a_{\ov{m-1}\,\ov{u-1}}$ (see
Figure 1). Similarly, $\bl_{\ov{m-1}}=\ov{u-1}\cdots \ov{u-1}$
($a_{\ov{m-1}\,\ov{u-1}}$ times). Repeating the above arguments, it
follows that for $0\leq k\leq u-1$,
\begin{equation*}
\bi_{\ov{u-k}}=\underbrace{\ov{m-k}\cdots \ov{m-k}}_{a_{\ov{m-k}\,\ov{u-k}}}\ ,
\ \ \ \ \bl_{\ov{m-k}}=\underbrace{\ov{u-k}\cdots \ov{u-k}}_{a_{\ov{m-k}\,\ov{u-k}}},
\end{equation*}
and
\begin{equation*}
\lambda_1\geq \cdots\geq \lambda_u,
\end{equation*}
where $\lambda_{k+1}=a_{\ov{m-k}\,\ov{u-k}}$ for $0\leq k<u$. \vskip
3mm

\textsc{Step 2.} Suppose that $a_{st}\neq 0$ for some $1\leq s\leq
n$ and $1\leq t\leq v$. Let $\ell$ be the maximum column index
($1\leq \ell\leq v$) such that $a_{s\ell}\neq 0$ for some $1\leq
s\leq n$.

First, we claim that $a_{st}=0$ for $s+t > \ell+1$.

Consider $a_{s\, \ell}$ for $1\leq s \leq n$. Suppose that $a_{s\,
\ell}\neq 0$ for some $s\geq 2$. Then we have ${e}_{s-1}A\neq 0$,
since $\epsilon^{(s-1)}(A)=(\cdots, +^{a_{s-1\,
\ell}},-^{a_{s\,\ell}})$ and there exists at least one $-$ in the
$(s-1)$-signature of $A$. This is a contradiction. So, we have $a_{1
\ell}\neq 0$ and $a_{s\ell}=0$ for $s\geq 2$. Next, consider $a_{s\,
\ell-1}$ for $1\leq s \leq n$. Suppose that $a_{s\, \ell-1}\neq 0$
for some $s\geq 3$. Then we also have ${e}_{s-1}A\neq 0$, since
$\epsilon^{(s-1)}(A)=(\cdots, +^{a_{s-1\,
\ell-1}},-^{a_{s\,\ell-1}})$, which is a contradiction. Similarly,
we can check that $a_{s t}=0$ for $1\leq t\leq \min{(\ell,n)}$ and
$s> \ell-t+1$.

Now, we claim that $a_{st}=a_{s't'}$ for $2\leq s+t=s'+t'\leq
\ell+1$.

Consider $a_{st}$ for $s+t=\ell+1$. Since ${e}_2 A=0$ and
${e}^*_{\ell-1}A=0$, we have $a_{2\, \ell-1}\leq a_{1\, \ell}$ and
$a_{2\, \ell-1}\geq a_{1\, \ell}$, respectively. Hence $a_{2\,
\ell-1}= a_{1\, \ell}$. Continuing this argument, it follows that
$a_{s\, \ell-s+1}=a_{1 \ell}$ for $1\leq s\leq \min{(\ell,n)}$.
Indeed, we have $\ell\leq n$. Otherwise, we have $a_{n\
\ell-n+1}\neq 0$, and $e_{\ell-n}^*A\neq 0$. Next, consider $a_{st}$
for $s+t=\ell$. Since  ${e}_1A=0$ and
$$\epsilon^{(1)}(A)=(\cdots, -^{a_{2\, \ell-2}},+^{a_{1\, \ell-1}},-^{a_{2\, \ell-1}},+^{a_{1\, \ell}})$$
(note that $a_{2\, \ell-1}=a_{1\, \ell}$), we have $a_{2\,
\ell-2}\leq a_{1\, \ell-1}$. Since ${e}_{\ell-2}^*A=0$ and
$$\epsilon^{(\ell-2)}(A)=(\cdots, -^{a_{1\, \ell-1}},+^{a_{2\, \ell-2}},-^{a_{2\, \ell-1}},+^{a_{3\, \ell-2}})$$
(note that $a_{3\, \ell-2}=a_{2\, \ell-1}=a_{1\, \ell}$), we have
$a_{2\, \ell-2}\geq a_{1\, \ell-1}$. Hence, $a_{2\, \ell-2}=
a_{1\, \ell-1}$. Similarly, we can check that $a_{1\,
\ell-1}=a_{2\, \ell-2}=a_{3\, \ell-3}=\cdots=a_{\ell-1\, 1}$.

Applying the above arguments successively, we conclude that $a_{s\,
t}=a_{s'\, t'}$ for $2\leq s+t=s'+t'\leq \ell+1$. We set
$\nu_k=a_{s\, k-s+1}$ for $1\leq k\leq \ell$.

\vskip 3mm

\textsc{Step 3}. Let $\ell$ be the maximum column index ($1\leq
\ell\leq v$) given in \textsc{Step 2}. We assume that $\ell=0$ if
$a_{st}= 0$ for all $1\leq s\leq n$ and $1\leq t\leq v$.

We claim that $a_{\ov{k} t}=1$ for $\ov{m-u}\leq \ov{k}\leq
\ov{1}$ and $1\leq t\leq \ell$. Let us use the induction on $t$.
If $\ell=0$, then it is clear. Suppose that $\ell>0$. Consider
$P_1=P(\bi_{\ov{u}}\cdots\bi_{\ov{1}}\bi_1)$. We have seen in
\textsc{Step 1} that $P(\bi_{\ov{u}}\cdots\bi_{\ov{1}})$ is an
$(m,n)$-hook semistandard tableau whose shape is a partition
$(\lambda_1,\cdots,\lambda_u)$ with the $k^{\rm th}$ row filled
with $\ov{m-k+1}$ ($1\leq k\leq u$). When we insert the word
$\bi_1$ into $P(\bi_{\ov{u}}\cdots\bi_{\ov{1}})$, all the letters
in $\bi_1$ are placed in the first column of $P_1$. If
$a_{\ov{k}\, 1}=0$ for some $\ov{m-u}\leq \ov{k}\leq \ov{1}$, then
there exists at least one letter in $\B_{m|n}^-$ placed in the
first $m$ rows of $P_1$ and hence $P(\bi)$. This contradicts the
fact that $P(\bi)$ is a highest weight tableau.

For $t<\ell$, suppose that $a_{\ov{k}t'}=1$ for $\ov{m-u}\leq
\ov{k}\leq \ov{1}$ and $1\leq t'\leq t$. Put
$P_t=P(\bi_{\ov{u}}\cdots\bi_{\ov{1}}\bi_1\cdots\bi_{t})$. Then,
each $k^{\rm th}$ row of $P_t$ ($1\leq k\leq m$) is filled with
$\ov{m-k+1}$. If we cut out the first $m$ rows of $P_t$, then the
remaining tableau consists of exactly $t$ columns. Moreover, if we
read its $k^{\rm th}$ column ($1\leq k\leq t$) from top to bottom,
then the associated word is given by
\begin{equation*}
k^{\nu_{k}+\cdots+\nu_{t}}(k+1)^{\nu_{t+1}}(k+2)^{\nu_{t+2}}\cdots(\ell-t+k)^{\nu_{\ell}}.
\end{equation*}
Since
\begin{equation*}
\bi_{t+1}=(\ov{m-u})^{a_{\ov{m-u}\, t+1}}\cdots \ov{1}^{a_{\ov{1}\,
t+1}}1^{\nu_{t+1}}\cdots (\ell-t)^{\nu_{\ell}},
\end{equation*}
it is not difficult to see that for $1\leq k\leq \ell-t$,
$\nu_{t+k}$ $(t+k)$'s are bumped out of the $t^{\rm th}$ column and
inserted into the $(t+1)^{\rm st}$ column, when we insert the word
$\bi_{t+1}$ into $P_t$ (note that $\nu_{\ell}>0$ and $\bi_{t+1}$ is
not an empty word). So, if $a_{\ov{k}\, t+1}=0$ for some
$\ov{m-u}\leq \ov{k}\leq \ov{1}$, then at least one letter in
$\B_{m|n}^-$ happens to be placed in the first $m$ rows of $P_{t+1}$
and hence $P(\bi)$, which is a contradiction. This completes the
induction.\vskip 3mm

\textsc{Step 4}. Consider $a_{\ov{1}\, t}$ for $\ell<t\leq v$.
Suppose that $a_{\ov{1}\, t}=1$ and $a_{\ov{1}\, t-1}=0$ for some
$\ell<t\leq v$. Then, we have ${e}_{t-1}A\neq 0$ since
$\epsilon^{(t-1)}(A)=(\cdots,+,-)$, which is a contradiction.

Next, consider $a_{\ov{2}\, t}$ for $\ell<t\leq v$. Suppose that
$a_{\ov{2}\, t}=1$ and $a_{\ov{2}\, t-1}=0$ for some $\ell+1<t\leq
v$. We assume that $t$ is the minimum index such that $a_{\ov{2}\,
t}=1$ and $a_{\ov{2}\, t-1}=0$. If $a_{\ov{1}\, t-1}=a_{\ov{1}\,
t}=0$ or $a_{\ov{1}\, t-1}=a_{\ov{1}\, t}=1$, then ${e}_{t-1}^*A\neq
0$. If $a_{\ov{1}\, t-1}=1$ and $a_{\ov{1}\, t}=0$, then $a_{\ov{1}
t'}=1$ for $t'< t$, and $a_{\ov{1} t''}=0$ for $t''\geq t$, which
implies that there exists a $\ov{1}$-good $-$ sign in $A$ with
respect to the $\frak{gl}_{m|n}$-crystal structure, and
${e}_{\ov{1}}A\neq 0$. So, in any case, we get a contradiction.
Moreover, we see that there is no $t$ ($\ell < t\leq v$) such that
$a_{\ov{2} t}=0$ and $a_{\ov{1} t}=1$, since ${e}_{\ov{1}}A= 0$.

Now, applying the above arguments successively to $a_{\ov{k}\, t}$
for $1< k\leq m-u$ and $\ell< t\leq v$, it follows that if
$a_{\ov{k} t}=0$ for $1< k\leq m-u$ and $\ell<t\leq v$, then
$a_{\ov{k}\, t+1}=a_{\ov{k-1}\, t}=0$. This leads to
\begin{equation*}
\lambda_{u+1}\geq\cdots\geq \lambda_{m},
\end{equation*}
where $\lambda_{u+k}=\sum_{1\leq t \leq v}a_{\ov{m-u-k+1}\, t}$
($1\leq k\leq m-u$). Also, we have $\lambda_u\geq \lambda_{u+1}$
since ${e}_{\ov{m-u}}A=0$.

Let $(\lambda_{m+1},\lambda_{m+2},\cdots,\lambda_r)$ be the
partition determined by \eqref{nu}, and put
\begin{equation*}
\lambda=(\lambda_1,\cdots,\lambda_m,\lambda_{m+1},\cdots,\lambda_r).
\end{equation*}
Note that $\lambda_{m+1}=\ell$, where $\ell$ is the maximum column
index given in \textsc{Step 2}, and $\lambda_{m}\geq \lambda_{m+1}$
by \textsc{Step 3}. So, $\lambda$ is a Young diagram, and
$\lambda\in \mathcal{P}_{m|n}\cap \mathcal{P}_{u|v}$ by
construction. Finally, we conclude that $A=A_{\lambda}$ given in
\eqref{Alambda}. \qed\vskip 3mm

Let $x=\{\,x_b\,|\,b\in\B_{m|n}\,\}$ and
$y=\{\,y_{b'}\,|\,b'\in\B_{u|v}\,\}$. The character of $\M$ is
given by
\begin{equation*}
{\rm ch}\M =\sum_{A\in\M}x^{{\rm wt}(A)}y^{{\rm wt}^*(A)}
=\frac{\prod_{|b|\neq |b'|}(1+x_{b}y_{b'})}
{\prod_{|b|=|b'|}(1-x_{b}y_{b'})},
\end{equation*}
where $b\in\B_{m|n}$ and $b'\in\B_{u|v}$. By Theorem \ref{highest
weight vector1}, we recover the {\it super Cauchy identity};
\begin{equation}\label{Cauchy1}
\frac{\prod_{|b|\neq |b'|}(1+x_{b}y_{b'})}
{\prod_{|b|=|b'|}(1-x_{b}y_{b'})}
=\sum_{\lambda\in\mathcal{P}_{m|n}\cap\mathcal{P}_{u|v}}hs_{\lambda}(x)hs_{\lambda}(y).
\end{equation}

\begin{rem}\label{remark on Kn map}{\rm
Let $\mathcal{A}_n=\{\,a_1<\cdots<a_n\,\}$ be the set of $n$ letters
with a linear ordering. For a partition $\lambda$ with length at
most $n$, a tableau $T$ obtained by filling $\lambda$ with the
letters in $\mathcal{A}_n$ is called a {\it semistandard tableau of
shape $\lambda$} if the entries in each row (resp. column) are
weakly (resp. strictly) increasing from left to right (resp. from
top to bottom). We denote by $\mathcal{SST}_n(\lambda)$ the set of
all semistandard tableau of shape $\lambda$ with entries in
$\mathcal{A}_n$. For example, $\B_{m|0}(\lambda)$
($\lambda\in\mathcal{P}_{m|0}$) may be identified with
$\mathcal{SST}_{m}(\lambda)$, and $\B_{0|n}(\lambda)$
($\lambda\in\mathcal{P}_{0|n}$) with $\mathcal{SST}_{n}(\lambda')$,
where $\lambda'$ is the conjugate of $\lambda$.

With this notation, we can recover several variations of the
original Knuth correspondence from Theorem \ref{highest weight
vector1} (cf.\cite{Fu}). If we put $n=v=0$ or $m=u=0$, then we have
two kinds of Knuth correspondence, where the one is given by the
column insertion of words and the other is given by the row
insertion of words. If we put $m=v=0$ or $n=u=0$, then we obtain the
dual Knuth correspondence. Similarly, one may also obtain other
variations from the decomposition given in next section (see Theorem
\ref{highest weight vector2}).

}
\end{rem}

\subsection{Diagonal action on symmetric matrices} Suppose that $m=u$
and $n=v$. Set
\begin{equation}
\mathfrak{M}=\{\,A\in\M\,|\,A=A^t\,\},
\end{equation}
the set of all symmetric matrices in $\M$. Let us consider the
diagonal action of $\frak{gl}_{m|n}\oplus\frak{gl}_{m|n}$ on
$\frak{M}$. That is, for $i\in I_{m|n}$ and $A\in\mathfrak{M}$, we
define
\begin{equation*}
\begin{split}
&{\bf e}_i A=e_ie_i^*A=e_i^*e_iA, \\
&{\bf f}_i A=f_if_i^*A=f_i^*f_iA.
\end{split}
\end{equation*}
Note that $P_1(A)=P_2(A)$ for $A\in \mathfrak{M}$ (see \eqref{pi}).
Put ${\rm wt}(A)={\rm wt}(P_1(A))={\rm wt}(P_2(A))$,
$\varepsilon_i(A)={\rm max}\{\,k\,|\,{\bf e}_i^k A\neq 0\,\}$, and
$\varphi_i(A)={\rm max}\{\,k\,|\,{\bf f}_i^kA\neq 0\,\}$  for $i\in
I_{m|n}$ and $A\in\mathfrak{M}$. Then we have

\begin{prop}\label{diagonal1} $\mathfrak{M}$ is a crystal graph for $\frak{gl}_{m|n}$, which decomposes as follows;
\begin{equation*}
\mathfrak{M} \simeq \bigoplus_{\lambda\in \mathcal{P}_{m|n}}
\B_{m|n}(\lambda).
\end{equation*}
\end{prop}
\pf If ${\bf e}_i A\neq 0$ or ${\bf f}_i A\neq 0$ for
$A\in\mathfrak{M}$ and $i\in I_{m|n}$, then
\begin{equation*}
\begin{split}
&({\bf e}_i A)^t=(e_ie_i^*A)^t=e_i^*e_iA^t={\bf e}_i A\in\mathfrak{M}, \\
&({\bf f}_i A)^t=(f_if_i^*A)^t=f_i^*f_iA^t={\bf f}_i
A\in\mathfrak{M}.
\end{split}
\end{equation*}
Hence, ${\bf e}_i,{\bf f}_i : \mathfrak{M} \rightarrow
\mathfrak{M}\cup \{0\}$ are well-defined operators for $i\in
I_{m|n}$.

For $A\in\mathfrak{M}$ and $i\in I_{m|n}$, we have $x_iA\neq 0$ if
and only if $x_i^*A\neq 0$ ($x=e,f$) since $A$ is symmetric. Hence,
by Lemma \ref{P-tableaux}, we have
\begin{equation*}
\begin{split}
&{\bf e}_iA \neq 0 \Leftrightarrow P_1(e_iA)\neq 0
\Leftrightarrow e_iP_1(A)\neq 0, \\
&{\bf f}_iA \neq 0 \Leftrightarrow P_1(f_iA)\neq 0 \Leftrightarrow
f_iP_1(A)\neq 0,
\end{split}
\end{equation*}
for $i\in I_{m|n}$. This implies that $\mathfrak{M}$ is a
$\frak{gl}_{m|n}$-crystal.

Next, consider the decomposition of $\mathfrak{M}$. For $A\in
\mathfrak{M}$, $A$ is $\frak{gl}_{m|n}$-equivalent to $P_1(A)$. So
each connected component in $\mathfrak{M}$ is generated by
$A_{\lambda}$ for some $\lambda\in\mathcal{P}_{m|n}$ by Theorem
\ref{highest weight vector1}. Since $A_{\lambda}\in \mathfrak{M}$
for $\lambda\in\mathcal{P}_{m|n}$, the set of all the highest weight
elements in $\mathfrak{M}$ is equal to $\M_{\rm h.w.}$. \qed\vskip
3mm

Now, let us show that there is an interesting relation between the
diagonal entries of a matrix in $\frak{M}$ and the shape of the
corresponding tableau (cf.\cite{Kn}), and hence obtain a family of
subcrystals of $\frak{M}$, which also have nice decompositions. For
$A=(a_{bb'})_{b,b'\in\B_{m|n}}\in\frak{M}$, let
\begin{equation}
\frak{o}(A)=\left|\,\{\,b\in\B_{m|n}^+\,|\,a_{bb}\equiv 1\hskip
-2mm\pmod{2}\,\}\,\right|+\sum_{b\in\B_{m|n}^-}a_{bb}.
\end{equation}

\begin{prop}\label{diagonal2}
Let $\frak{M}_k=\{\,A\in\frak{M}\,|\,\frak{o}(A)=k\,\}$ for $k\geq
0$. Then $\frak{M}_k$ is a subcrystal of $\frak{M}$, and decomposes
as follows;
$$\frak{M}_k\simeq
\bigoplus_{\substack{\lambda \in\mathcal{P}_{m|n} \\
o(\lambda)=k}}\B_{m|n}(\lambda),$$ where $o(\lambda)$ is the number
of odd parts in $\lambda$.
\end{prop}
\pf Fix $k\geq 0$. First, we will check that $\mathfrak{M}_k$
together with $0$ is stable under ${\bf e}_i, {\bf f}_i$ $(i\in
I_{m|n})$, which implies that $\frak{M}_k$ is a subcrystal of
$\frak{M}$.

Given $A=(a_{bb'})\in\frak{M}_k$ and $i\in I_{m|n}$, suppose that
${\bf f}_iA\neq 0$. We assume that the diagonal entries of $A$ are
changed under ${\bf f}_i$, equivalently under $f_i$ or $f^*_i$.
Otherwise, it is clear that ${\bf f}_iA\in\frak{M}_k$. For
convenience,
we write $$A[i]=\left(%
\begin{array}{cc}
  a_{bb} & a_{bb'} \\
  a_{b'b} & a_{b'b'} \\
\end{array}%
\right),$$ where $b<b'$ are the indices such that $\langle
h_i,\epsilon_b\rangle, \langle h_i,\epsilon_{b'}\rangle\neq
0$.\vskip 3mm

{\textsc Case 1}. $i = \ov{k}$, ($1\leq k \leq m-1$). Consider
$$A[i]=\left(%
\begin{array}{cc}
  a_{\ov{k+1}\,\ov{k+1}} & a_{\ov{k+1}\,\ov{k}} \\
  a_{\ov{k}\,\ov{k+1}} & a_{\ov{k}\,\ov{k}} \\
\end{array}%
\right)=\left(%
\begin{array}{cc}
  a & b \\
  b & c \\
\end{array}%
\right).$$ Suppose that
$(f_{i}A)[i]=\left(%
\begin{array}{cc}
  a-1 & b \\
  b+1 & c \\
\end{array}%
\right).$ This implies that with respect to $f_i$, we have
$\epsilon^{(i)}(A)=(\cdots,-^b,+^a,-^c,+^b ,\cdots)$, where $a>c$
and the $i$-good $+$ sign of $A$ appears in $+^a$ of
$\epsilon^{(i)}(A)$. Note that
$\epsilon^{(i)}(f_iA)=(\cdots,-^b,+^{a-1},-^{c},+^{b+1} ,\cdots)$
with respect to  $f_i^*$.

If $a\equiv c \pmod{2}$, then we still have $a-1>c$, and the
$i$-good $+$ sign of $f_iA$ with respect to $f_i^*$ appears in
$+^{a-1}$, and we have
$$({\bf f}_iA)[i]=(f^*_if_{i}A)[i]=\left(%
\begin{array}{cc}
  a-2 & b+1 \\
  b+1 & c \\
\end{array}%
\right).$$ If $a\not\equiv c \pmod{2}$, then we have
$$({\bf f}_iA)[i]=(f^*_if_{i}A)[i]=\left(%
\begin{array}{cc}
  a-2 & b+1 \\
  b+1 & c \\
\end{array}%
\right) \ \ \text{or} \ \ \left(
\begin{array}{cc}
  a-1 & b \\
  b & c+1 \\
\end{array}
\right).$$ In any case, we have $\frak{o}({\bf f}_iA)=k$ and ${\bf
f}_iA\in\frak{M}_k$.

Next, suppose that
$(f_{i}A)[i]=\left(%
\begin{array}{cc}
  a & b-1 \\
  b & c+1 \\
\end{array}%
\right)$. Then we must have $({\bf f}_iA)[i]=\left(%
\begin{array}{cc}
  a & b-1 \\
  b-1 & c+2 \\
\end{array}%
\right)$, and $\frak{o}({\bf f}_iA)=k$.\vskip 3mm

{\textsc Case 2}. $i = k$, ($1\leq k \leq n-1$). Consider
$$A[i]=\left(%
\begin{array}{cc}
  a_{k\, k} & a_{k\,k+1} \\
  a_{k+1\, k} & a_{k+1\,k+1} \\
\end{array}%
\right)=\left(%
\begin{array}{cc}
  a & b \\
  b & c \\
\end{array}%
\right).$$ Suppose that $a>0$ and
$(f_{i}A)[i]=\left(%
\begin{array}{cc}
  a-1 & b \\
  b+1 & c \\
\end{array}%
\right).$ Then with respect to $f_i$, we have
$\epsilon^{(i)}(A)=(\cdots,+^a,-^b,+^b,-^c,\cdots)$, where the
$i$-good $+$ sign of $A$ appears in $+^a$. Since
$\epsilon^{(i)}(f_iA)=(\cdots,+^{a-1},-^{b},+^{b+1},-^{c} ,\cdots)$
with respect to $f_i^*$, the $i$-good $+$ sign of $f_iA$ appears in
$+^{b+1}$, and we have
$$({\bf f}_iA)[i]=(f^*_if_{i}A)[i]=\left(%
\begin{array}{cc}
  a-1 & b \\
  b & c+1 \\
\end{array}%
\right),$$ which implies that $\frak{o}({\bf f}_iA)=k$.

Next, suppose that $b>0$ and
$(f_{i}A)[i]=\left(%
\begin{array}{cc}
  a & b-1 \\
  b & c+1 \\
\end{array}%
\right)$. But, this can't happen since
$\epsilon^{(i)}(A)=(\cdots,+^a,-^b,+^b,-^c,\cdots)$ with respect to
$f_i$, and the pair $(-^b,+^b)$ cancels out. \vskip 3mm

\textsc{Case 3}. $i = 0$. Consider
$$A[0]=\left(%
\begin{array}{cc}
  a_{\ov{1}\, \ov{1}} & a_{\ov{1}\,1} \\
  a_{1\, \ov{1}} & a_{1\,1} \\
\end{array}%
\right)=\left(%
\begin{array}{cc}
  a & b \\
  b & c \\
\end{array}%
\right).$$ Then we have
$$({\bf f}_0A)[0]=(f^*_0f_{0}A)[0]=\left(%
\begin{array}{cc}
  a-2 & b+1 \\
  b+1 & c \\
\end{array}%
\right) \ \ \text{or} \ \ \left(
\begin{array}{cc}
  a-1 & b \\
  b & c+1 \\
\end{array}
\right),$$ (it can't happen that $({\bf f}_0A)[0]=\left(%
\begin{array}{cc}
  a & b-1 \\
  b-1 & c+2 \\
\end{array}%
\right)$). In any case, we have $\frak{o}({\bf f}_iA)=k$.

Similarly, we can check that ${\bf e}_i\frak{M}_k\subset
\frak{M}_k\cup\{0\}$ for $i\in I_{m|n}$. Therefore, $\frak{M}_k$ is
a crystal graph for $\frak{gl}_{m|n}$.\vskip 3mm

Next, we observe that for $\lambda=(\lambda_k)_{k\geq
1}\in\mathcal{P}_{m|n}$,
$$\frak{o}(A_{\lambda})=\left|\,\{\,i\,|\,1\leq i\leq m, \, \lambda_i \ \text{is odd}\,\}\,\right|+\sum_{k\geq 0}\nu_{2k+1},$$ (see
\eqref{Alambda}). Since $\nu_i$ is the number of occurrences of $i$
in $(\lambda_{m+1},\lambda_{m+2},\cdots)$ (see \eqref{nu}), it
follows that $\frak{o}(A_{\lambda})$ is the number of odd parts in
$\lambda$, say $o(\lambda)$, and hence
$A_{\lambda}\in\mathfrak{M}_k$ if and only if $o(\lambda)=k$.
\qed\vskip 3mm

\begin{cor}\label{diagonal3} Under the above hypothesis, we have
$\frak{M}=\bigoplus_{k\geq 0}\frak{M}_k$, and in particular,
$$\frak{M}_0\simeq
\bigoplus_{\substack{\lambda \in\mathcal{P}_{m|n} \\
\lambda\,:\,\text{even}}}\B_{m|n}(\lambda).$$
\end{cor}
\pf It follows from the fact that $\frak{o}(A_{\lambda})=0$ if and
only if $\lambda$ is even (that is, each part of $\lambda$ is even).
\qed

\begin{rem}{\rm
(1) A special case of Proposition \ref{diagonal2} was first observed
in \cite{Kn}. Let us give a brief explanation. Put $m=0$ in
Corollary \ref{diagonal3}. We identify $\B_{0|n}$ with
$\mathcal{A}_n$ (see Remark \ref{remark on Kn map}), and
$\B_{0|n}(\lambda)$ with $\mathcal{SST}_n(\lambda')$.  Hence, the
set of all $n\times n$ symmetric matrices of non-negative integers
with ${\rm tr}(A)=k$, is in one-to-one correspondence with
$\bigsqcup_{o(\lambda')=k}\mathcal{SST}_n(\lambda)$ (Theorem 4
\cite{Kn}) where the union is given over all partitions with the
number of odd columns $k$.

(2) If we consider the characters of the decompositions in
Proposition \ref{diagonal1} and Corollary \ref{diagonal3}, then we
obtain the following identities (cf.\cite{Mac});
\begin{equation*}
\dfrac{\prod_{b<b',\,
|b|\neq|b'|}(1+x_bx_{b'})}{\prod(1-x_b)\prod_{b<b',\,
|b|=|b'|}(1-x_bx_{b'})}
=\sum_{\lambda\in\mathcal{P}_{m|n}}hs_{\lambda}(x),
\end{equation*}
where $b,b'\in\B_{m|n}$,
\begin{equation*}
\dfrac{\prod_{b<b',\, |b|\neq|b'|}(1+x_bx_{b'})}{\prod_{|b|=0
}(1-x_{b}^2)\prod_{b< b',\, |b|=|b'|}(1-x_bx_{b'})}
=\sum_{\substack{\lambda\in\mathcal{P}_{m|n} \\ \lambda
\,:\,\text{\rm even}}}hs_{\lambda}(x),
\end{equation*}
where   $b,b'\in\B_{m|n}$.

}
\end{rem}

\section{Dual construction}
In this section, we discuss a bicrystal graph  associated to the
super exterior algebra $\Lambda(\mathbb{C}^{m|n}\otimes
\mathbb{C}^{u|v})$, and its explicit decomposition.

Suppose that  $m,n,u,v$ are non-negative integers such that $m+n,
u+v>0$. We set
\begin{equation}
\begin{split}
\M^{\sharp}_{m|n,u|v}=\{\,A=&(a_{bb'})_{b\in \B_{m|n}, b'\in
\B_{u|v}}\,|\,  \\ &\text{(1) $a_{bb'}\in\mathbb{Z}_{\geq 0}$,\ \
(2) $a_{bb'}\leq 1$ if $|b|=|b'|$}\, \}.
\end{split}
\end{equation}
For convenience, we write $\M^{\sharp}=\M^{\sharp}_{m|n,u|v}$.

As in the case of $\M$, we identify a matrix in $\M^{\sharp}$ with
a biword given by reading the row and column indices of non-zero
entries of the matrix with respect to a linear ordering. First, we
let
\begin{equation}
\begin{split}
\Omega^{\sharp}_{m|n,u|v}=\Omega^{\sharp}=\{\, &(\bi,\bj)\in
\mathcal{W}_{m|n}\times
\mathcal{W}_{u|v}\, | \, \\
& \text{(1) $\bi=i_1\cdots i_r$ and $\bj=j_1 \cdots j_r$ for some
$r\geq 0$,} \\
& \text{(2) $(i_1,j_1)\preceq \cdots \preceq (i_r,j_r)$,} \\
& \text{(3) $|i_k|= |j_k|$ implies $(i_k,j_k)\neq (i_{k\pm1},j_{k\pm
1})$,} \,\},
\end{split}
\end{equation}
where for $(i,j)$ and $(k,l)\in \B_{m|n}\times \B_{u|v}$, the
linear ordering $\prec$ is given by
\begin{equation}\label{partial order}
(i,j)\prec (k,l) \ \ \ \ \Leftrightarrow \ \ \ \
\begin{cases}
(j<l) & \text{or}, \\
(j=l\in \B_{u|v}^+,\ \text{and} \ i<k) & \text{or}, \\
(j=l\in \B_{u|v}^-,\ \text{and} \ i>k) &.
\end{cases}
\end{equation}
We define $e_i,f_i : \Omega^{\sharp} \rightarrow
\Omega^{\sharp}\cup \{0\}$ ($i\in I_{m|n}$) by
\begin{equation*}
e_i (\bi,\bj)=(\re_i\bi,\bj), \ \ \ \ \ f_i
(\bi,\bj)=(\rf_i\bi,\bj),
\end{equation*}
for $(\bi,\bj)\in \Omega^{\sharp}$. Set ${\rm wt}(\bi,\bj)={\rm
wt}(\bi)$, $\varepsilon_i(\bi,\bj)=\varepsilon_i(\bi)$ and
$\varphi_i(\bi,\bj)=\varphi_i(\bi)$ ($i\in I_{m|n}$). Then it is
easy to see that $\Omega^{\sharp}$ is a crystal graph for
$\frak{gl}_{m|n}$ (cf. Lemma \ref{crystal}).

For $(\bi,\bj)\in \Omega^{\sharp}_{m|n,u|v}$, we define
$A(\bi,\bj)=(a_{bb'})$ to be a matrix in $\M^{\sharp}$, where
$a_{bb'}$ is the number of $k$'s such that $(i_k,j_k)=(b,b')$ for
$b\in \B_{m|n}$ and $b'\in \B_{u|v}$. Then, the map
$(\bi,\bj)\mapsto A(\bi,\bj)$ is a bijection between
$\Omega^{\sharp}_{m|n,u|v}$ and $\M^{\sharp}$. Hence, $\M^{\sharp}$
is a crystal graph for $\frak{gl}_{m|n}$ with this identification.

Next, we introduce $(\Omega^{\sharp}_{m|n,u|v})^*$  to define a
$\frak{gl}_{u|v}$-crystal structure on $\M^{\sharp}$. Recall that in
Section 3, a biword in $\Omega^*$ was obtained by reading the row
and column indices of the non-zero entries in the transpose of a
given matrix with respect to the same lexicographic ordering used in
$\Omega$. But in the case of $\M^{\sharp}$, we need another linear
ordering. That is, we set
\begin{equation}
\begin{split}
(\Omega^{\sharp}_{m|n,u|v})^*=(\Omega^{\sharp})^*=\{\, &(\bk,\bl)\in
\mathcal{W}_{m|n}\times
\mathcal{W}_{u|v}\, | \, \\
& \text{(1) $\bk=k_1 \cdots k_r$ and $\bl=l_1 \cdots l_r$ for some
$r\geq 0$,} \\
& \text{(2) $(k_1,l_1)\preceq' \cdots \preceq' (k_r,l_r)$,} \\
& \text{(3) $|k_t|= |l_t|$ implies $(k_t,l_t)\neq (k_{t\pm1},l_{t\pm
1})$,} \,\},
\end{split}
\end{equation}
where for $(i,j)$ and $(k,l)\in \B_{m|n}\times \B_{u|v}$, the linear
ordering $\prec'$ is given by
\begin{equation}\label{partial order}
(i,j)\prec' (k,l) \ \ \ \ \Leftrightarrow \ \ \ \
\begin{cases}
(i>k) & \text{or}, \\
(i=k\in \B_{u|v}^+,\ \text{and} \ j<l) & \text{or}, \\
(i=k\in \B_{u|v}^-,\ \text{and} \ j>l) &.
\end{cases}
\end{equation}
Clearly, we have a bijection $(\bk,\bl)\mapsto A(\bk,\bl)$ from
$(\Omega^{\sharp})^*$ to $\M^{\sharp}$.

Similarly, we define  $e_j^*,f_j^* : (\Omega^{\sharp})^*
\rightarrow (\Omega^{\sharp})^*\cup \{0\}$ ($j\in I_{u|v}$) by
\begin{equation*}
e_j^* (\bk,\bl)=(\bk,\re_j\bl), \ \ \ \ \ f_j^*
(\bk,\bl)=(\bk,\rf_j\bl),
\end{equation*}
for $(\bk,\bl)\in (\Omega^{\sharp})^*$. Set ${\rm
wt}^*(\bk,\bl)={\rm wt}(\bl)$,
$\varepsilon^*_j(\bk,\bl)=\varepsilon_j(\bl)$ and
$\varphi^*_j(\bk,\bl)=\varphi_j(\bl)$ ($j\in I_{u|v}$). Then
$(\Omega^{\sharp})^*$ is a crystal graph for $\frak{gl}_{u|v}$
(cf. Lemma \ref{crystal}), and hence so is $\M^{\sharp}$.

\begin{ex}\label{example2}{\rm
Suppose that $m|n=u|v=2|2$ and
\begin{equation*}
A=\left(\begin{tabular}{cc|cc}
 1 & 1 & 0 & 0   \\
 0 & 0 & 2 & 1   \\ \hline
 0 & 1 & 1 & 0   \\
 2 & 0 & 0 & 0   \\
\end{tabular}\right )\in\M^{\sharp}.
\end{equation*}
Then $A=A(\bi,\bj)=A(\bk,\bl)$ for $(\bi,\bj)\in\Omega^{\sharp}$
and $(\bk,\bl)\in(\Omega^{\sharp})^*$, where
\begin{equation*}
\begin{split}
&\bi=\ov{2} \ 2 \ 2 \ \ov{2} \ 1 \ 1 \ \ov{1} \ \ov{1} \ \ov{1}, \ \ \ \
\bk=2 \ 2 \ 1 \ 1 \ \ov{1} \ \ov{1} \ \ov{1} \ \ov{2} \ \ov{2},  \\
&\bj=\ov{2} \ \ov{2} \ \ov{2} \ \ov{1} \ \ov{1} \ 1 \ 1 \ 1 \ 2, \ \ \ \
\bl\,\, =\ov{2} \ \ov{2} \ 1 \ \ov{1} \ 1 \ 1 \ 2 \ \ov{2} \ \ov{1}.
\end{split}
\end{equation*}
}
\end{ex}

For $A\in\M^{\sharp}$, we have $A=A(\bi,\bj)=A(\bk,\bl)$ for unique
$(\bi,\bj)\in \Omega^{\sharp}$ and
$(\bk,\bl)\in(\Omega^{\sharp})^*$. Then we define
\begin{equation}\label{pisharp}
{\pi}^{\sharp}(A)=(P_1^{\sharp}(A),P_2^{\sharp}(A))=(P(\bi),P(\bl)).
\end{equation}
By definition, we have $A\simeq_{\frak{gl}_{m|n}}P_1^{\sharp}(A)$,
and $A\simeq _{\frak{gl}_{u|v}}P_2^{\sharp}(A)$.

Now, we have following analogue of Lemma \ref{bicrystal} and
Proposition \ref{g-crystal};
\begin{prop}\label{bicrystal'}
$\M^{\sharp}$ is a $(\frak{gl}_{m|n},\frak{gl}_{u|v})$-bicrystal,
and for each connected component $\mathbf{C}$ in $\M^{\sharp}$,
$\pi^{\sharp}$ gives the following isomorphism of
$(\frak{gl}_{m|n},\frak{gl}_{u|v})$-bicrystals;
$$\pi^{\sharp} : \mathbf{C} \longrightarrow \B_{m|n}(\lambda)\times \B_{u|v}(\mu)$$
for some $\lambda\in\mathcal{P}_{m|n}$ and
$\mu\in\mathcal{P}_{u|v}$.
\end{prop}
\pf For $A\in \M^{\sharp}$, let $A^{\rho}$ be the matrix given as a
clock-wise rotation of $A$ by $90^{\circ}$ (see Figure 2). The map
$A\mapsto A^{\rho}$ gives a one-to-one correspondence from
$\M^{\sharp}_{m|n,u|v}$ to $\M_{u|v,n|m}$. Moreover, we have
\begin{equation*}
\begin{split}
&(e_{\overline{k}}A)^{\rho}=f^*_{k}A^{\rho}, \ \
(f_{\overline{k}}A)^{\rho}=e^*_{k}A^{\rho}, \ \ (1\leq k\leq m-1),
\\
&(e_{l}A)^{\rho}=f^*_{\overline{l}}A^{\rho}, \ \
(f_{l}A)^{\rho}=e^*_{\overline{l}}A^{\rho},\ \ (1\leq l\leq n-1),\\
&(e_{0}A)^{\rho}=f^*_{0}A^{\rho}, \ \
(f_{0}A)^{\rho}=e^*_{0}A^{\rho},
\end{split}
\end{equation*}
and $(x^*_jA)^{\rho}=x_jA^{\rho}$ for $j\in I_{u|v}$ and $x=e,f$. By
Lemma \ref{bicrystal}, $e_i,f_i$ ($i\in I_{m|n}$) commute with
$e_j^*,f_j^*$ ($j\in I_{u|v}$) on $\M^{\sharp}$. Hence,
$\M^{\sharp}$ is a $(\frak{gl}_{m|n},\frak{gl}_{u|v})$-bicrystal.

Applying the same arguments in Lemma \ref{connected component},
\ref{P-tableaux} and Proposition \ref{g-crystal}, we conclude that
each connected component of $\M^{\sharp}$ is isomorphic to
$\B_{m|n}(\lambda)\times \B_{u|v}(\mu)$ for some
$\lambda\in\mathcal{P}_{m|n}$ and $\mu\in\mathcal{P}_{u|v}$. \qed

\begin{figure}
$\left(\begin{tabular}{ccc|ccc}
0 & 1 & 1 & 0 & 0 & 2   \\
1 & 0 & 0 & 2 & 1 & 0 \\ \hline
3 & 0 & 1 & 1 & 0 & 1 \\
0 & 2 & 0 & 0 & 0 & 0 \\
\end{tabular}\right )^{\rho}=
\left(\begin{tabular}{cc|cc}
 0 & 3 & 1 & 0   \\
 2 & 0 & 0 & 1   \\
 0 & 1 & 0 & 1   \\ \hline
 0 & 1 & 2 & 0   \\
 0 & 0 & 1 & 0   \\
 0 & 1 & 0 & 2
\end{tabular}\right )$
\caption{A clock-wise  rotation of $A$ by $90^{\circ}$}
\end{figure}

\begin{ex}{\rm
Let $A$ be the matrix given in Example \ref{example2}. Then we have
\begin{equation*}
{\pi}^{\sharp}(A)=\left(\begin{array}{ccccc}
\ov{2} & \ov{2} & \ov{1} & 1 & 2 \\
\ov{1} & \ov{1} &  &  &    \\
1 & 2 &  &  &
\end{array},
\begin{array}{ccc}
\ov{2} & \ov{2} & \ov{2}  \\
\ov{1} & \ov{1} & \ov{1}  \\
1 &  &    \\
1 &  &    \\
2 &  &
\end{array}\right ).
\end{equation*}

}
\end{ex}

Set
\begin{equation}
\M^{\sharp}_{\rm
h.w.}=\{\,A\in\M^{\sharp}\,|\,{\pi}^{\sharp}(A)=(H^{\lambda}_{m|n},H^{\mu}_{u|v})\
\text{for some
$(\lambda,\mu)\in\mathcal{P}_{m|n}\times\mathcal{P}_{u|v}$}\,\},
\end{equation}
which is the set of all the  highest weight elements in
$\M^{\sharp}$.

Suppose that $\lambda=(\lambda_k)_{k\geq
1}\in\mathcal{P}_{m|n}\cap\mathcal{P}_{v|u}$ is given, and let
$\lambda'=(\lambda'_k)_{k\geq 1}$ be the conjugate of $\lambda$. We
define ${A}^{\sharp}_{\lambda}=({a}^{\sharp}_{bb'})\in\M^{\sharp}$
by
\begin{itemize}
\item[(1)] for $0\leq k< m$, $0\leq l<u$,
\begin{equation*}
{a}^{\sharp}_{\ov{m-k}\, \ov{u-l}}=
\begin{cases}
1, & \text{if $\lambda_{k+1}\geq l+1$}, \\
0, & \text{otherwise}.
\end{cases}
\end{equation*}

\item[(2)] for $0\leq k< m$, $1\leq t\leq v$,
\begin{equation*}
{a}^{\sharp}_{\ov{m-k}\, t}=
\begin{cases}
\mu_t, & \text{if $k+1=t$}, \\
0, & \text{otherwise},
\end{cases}
\end{equation*}
where $\mu_t=\max{(\lambda_t-u,0)}$ for $t\geq 1$.

\item[(3)] for $1\leq s\leq n$, $0\leq l< u$,
\begin{equation*}
{a}^{\sharp}_{s\, \ov{u-l}}=
\begin{cases}
\nu_{s+l}, & \text{if $1\leq s+l \leq\ell$}, \\
0, & \text{otherwise},
\end{cases}
\end{equation*}
where $\nu_1,\cdots,\nu_{\ell}$ are given in \eqref{nu}.

\item[(4)]  for $1\leq s\leq n$, $1\leq t\leq v$,
${a}^{\sharp}_{st}=0$.
\end{itemize}

Note that ${A}^{\sharp}_{\lambda}$ is of the following form;
\begin{equation}\label{Alambda2}
{A}^{\sharp}_{\lambda}=\left(%
\begin{array}{cccccccc|ccccc}
   1 & 1 &  & \cdots& & & & 1 &  \mu_1 & & & &  \\
   1 & 1 &  &\cdots & & & & 1 & & \mu_2 & & &  \\
   \vdots & \vdots & & & & & & \vdots & & & \mu_3 & &  \\
   1 & 1 & & \cdots& & & & 1 & & & & \ddots &  \\
   \vdots& & & & & & & & & & & &   \\
   1& 1 & \cdots &  & 1 & & & & & & & &   \\ \hline
   \nu_{1} & \nu_{2} & \cdots & \nu_{\ell} & & & &  & & & & &    \\
   \nu_{2} & \nu_{3} & \cdots & & & & & & & & & &    \\
   \vdots & & & & & & & & & & & &    \\
   \nu_{\ell} & & & & & & & & & & & &    \\
   & & & & & & & & & & & &
\end{array}%
\right).
\end{equation}
By Schensted's algorithm, we can check that
\begin{equation*}
{\pi}^{\sharp}({A}^{\sharp}_{\lambda})=(H_{m|n}^{\lambda},H^{\lambda'}_{u|v}).
\end{equation*}
\begin{ex}{\rm
Suppose that $m|n=3|2$ and $u|v=4|3$. Let $\lambda=(7,6,2,2,1,1)$.
Then
\begin{equation*}
{A}^{\sharp}_{\lambda}=\left(%
\begin{array}{cccc|ccc}
  1 & 1 & 1 & 1 & 3 & 0 & 0 \\
  1 & 1 & 1 & 1 & 0 & 2 & 0 \\
  1 & 1 & 0 & 0 & 0 & 0 & 0 \\ \hline
  2 & 1 & 0 & 0 & 0 & 0 & 0 \\
  1 & 0 & 0 & 0 & 0 & 0 & 0
\end{array}%
\right),
\end{equation*}
and
\begin{equation*}
{\pi}^{\sharp}({A}^{\sharp}_{\lambda})=\left(
\begin{array}{ccccccc}
  \ov{3}  & \ov{3}  & \ov{3}  & \ov{3}  & \ov{3}  & \ov{3} & \ov{3} \\
  \ov{2}  & \ov{2}  & \ov{2}  & \ov{2}  & \ov{2}  & \ov{2} & \\
  \ov{1}  & \ov{1}  &   &   &   &  & \\
  1  & 2  &   &   &   &  & \\
  1  &   &   &   &   &  & \\
  1  &   &   &   &   &  &
\end{array},
\begin{array}{ccccccc}
  \ov{4}  & \ov{4}  & \ov{4}  & \ov{4}  & \ov{4}  & \ov{4} &  \\
  \ov{3}  & \ov{3}  & \ov{3}  & \ov{3}  &  &  & \\
  \ov{2}  & \ov{2}  &    &  &  &  & \\
  \ov{1}  & \ov{1}  &   &   &   &  & \\
  1  & 2  &   &   &   &  & \\
  1  & 2  &   &   &   &  & \\
  1  &   &   &   &   &  &
\end{array} \right ).
\end{equation*}

}
\end{ex}

Now, we can characterize all the highest weight elements in
$\M^{\sharp}$. The proof is almost the same as in Theorem
\ref{highest weight vector1}.
\begin{thm}\label{highest weight vector2} We have
\begin{equation*}
\M^{\sharp}_{\rm h.w.}=\{\,{A}^{\sharp}_{\lambda}\,|\,\lambda\in
\mathcal{P}_{m|n}\cap\mathcal{P}_{v|u}\,\},
\end{equation*}
and hence the following isomorphism of
$(\frak{gl}_{m|n},\frak{gl}_{u|v})$-bicrystals;
\begin{equation*}
\pi^{\sharp} : \M^{\sharp}\longrightarrow \bigoplus_{\lambda\in
\mathcal{P}_{m|n}\cap\mathcal{P}_{v|u}} \B_{m|n}(\lambda)\times
\B_{u|v}(\lambda').
\end{equation*}\qed
\end{thm}

In terms of characters, we also recover the {\it dual Cauchy
identity of hook Schur functions};
\begin{equation}
\frac{\prod_{|b|=|b'|}(1+x_{b}y_{b'})}{\prod_{|b|\neq|b'|}(1-x_{b}y_{b'})}
=\sum_{\lambda\in\mathcal{P}_{m|n}\cap\mathcal{P}_{v|u}}hs_{\lambda}(x)hs_{\lambda'}(y),
\end{equation}
where $b\in\B_{m|n}$ and $b'\in\B_{u|v}$.

\section{Semi-infinite construction }
Let $\frak{g}$ be a contragredinet Lie superalgebra of infinite rank
whose associated Dynkin diagram is given by
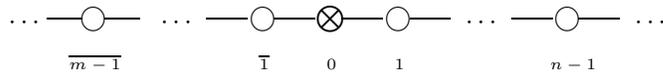
\begin{center}
\setlength{\unitlength}{0.45cm}
\begin{picture}(15,2)(0,0)
\put(-.85,.5){\line(1,0){1}}
\put(0,0){\makebox(1,1){$\bigcirc$}}\put(.85,.5){\line(1,0){1}}
\put(5,0){\makebox(1,1){$\bigcirc$}}\put(3.85,.5){\line(1,0){1.2}}\put(5.85,.5){\line(1,0){1.2}}
\put(7,0){\makebox(1,1){$\bigotimes$}}\put(7.85,.5){\line(1,0){1.2}}
\put(9,0){\makebox(1,1){$\bigcirc$}}\put(9.85,.5){\line(1,0){1.2}}
\put(14,0){\makebox(1,1){$\bigcirc$}}\put(12.85,.5){\line(1,0){1.2}}\put(14.85,.5){\line(1,0){1.2}}
\put(2.5,.2){$\cdots$}\put(11.5,.2){$\cdots$}
\put(-2,.2){$\cdots$}\put(16.5,.2){$\cdots$}

\put(-.2,-1){\tiny $\ov{m-1}$}\put(5.4,-1){\tiny
$\ov{1}$}\put(7.4,-1){\tiny $0$}\put(9.4,-1){\tiny $1$}\put(14
,-1){\tiny $n-1$}
\end{picture}\ \ \ \ \ \ \ \ \ \ \
\end{center} \vskip .5cm
(see \cite{Kac}). Then for all $m,n\geq 1$, there is a natural
embedding of $\frak{gl}_{m|n}$ into $\frak{g}$. Note that $\frak{g}$
is not equal to $\frak{gl}_{\infty|\infty}$ in the sense of
\cite{KacL}, but a proper subalgebra of it. In this section, we
study crystal graphs for $\frak{g}$ which are generated by highest
weight elements. Since a tensor power of the $\frak{g}$-crystal
associated to the natural representation does not have a highest
weight element, we introduce a $\frak{g}$-crystal $\mathscr{F}$
consisting of semi-infinite words, which is analogous to the crystal
graph for $\frak{gl}_{\infty}$ associated to a level one fermionic
Fock space representation (cf.\cite{Kac90}). We will show that each
connected component of a tensor power of $\mathscr{F}$ can be
realized as the set of semi-infinite semistandard tableaux, which is
generated by a highest weight vector. Then by using the methods
developed in the previous sections, we give an explicit
decomposition of $\mathscr{F}^{\otimes u}$ ($u\geq 2$) as a
$(\frak{g},\frak{gl}_{u})$-bicrystal.

\subsection{Crystal graphs of semi-infinite words}
We may naturally define a crystal graph for $\frak{g}$ by taking $m$
and $n$ to infinity in Definition \ref{crystal graph}. Set
\begin{equation*}
\B=\{\,\cdots<\ov{m}<\cdots<\ov{1}<1<\cdots<{n}<\cdots\,\},
\end{equation*}
and $\B^{+}$ (resp. $\B^{-}$) denotes the set of elements with
degree $0$ (resp. $1$) in $\B$. The index set is given by
\begin{equation*}
 I=\{\,\cdots,\ov{m},\cdots,\ov{1},0,1,\cdots,{n},\cdots\,\}.
\end{equation*}
Then the simple root $\alpha_i$  and the simple coroot $h_i$ ($i\in
I$) are defined in the same way. But, instead of
$\bigoplus_{b\in\B}\mathbb{Z}\epsilon_b$, we use
$$P=\mathbb{Z}\Lambda\oplus\bigoplus_{b\in\B}\mathbb{Z}\epsilon_b$$ as
the weight lattice of $\frak{g}$, where $\Lambda=\Lambda_0\in
(\bigoplus_{i\in I}\mathbb{Z}h_i)^*$ is the fundamental weight such
that $\langle h_i,\Lambda\rangle=\delta_{0 i}$ for $i\in I$. Note
that $\B$ is the crystal graph for $\frak{g}$ associated to the
natural representation of $\frak{g}$.

Now, we define $\mathscr{F}$ to be the set of {\it semi-infinite}
words $w=\cdots w_3w_2w_1$ with letters in $\B$ such that
\begin{itemize}
\item[(1)] $w_{i+1}\leq w_i$ for $i\geq 1$,

\item[(2)] $w_{i+1}=w_i$ implies $|w_i|=1$,

\item[(3)] there exists $c\in\mathbb{Z}$ such that $w_i=\ov{i+c}$ for all $i\gg 1$.
\end{itemize}
We call $c$ in (3) the {\it charge of $w$}. For example, the charge
of $w=\cdots\ov{6}\,\ov{5}\,\ov{4}\,\ov{3}\,\ov{1}\,2\, 3\, 4\,
4\in\mathscr{F}$ is $-3$. For each $i\in I$, we define the Kashiwara
operators
\begin{equation*}
\re_i,\rf_i : \mathscr{F} \longrightarrow  \mathscr{F}\cup\{0\}
\end{equation*}
as in the case of $\mathcal{W}_{m|n}$ (see Section 2.2). They are
well-defined since for $w=\cdots w_3w_2w_1\in\mathscr{F}$ and $i\in
I$,
$$\epsilon^{(i)}(w)=(\cdots,\epsilon^{(i)}(w_3),\epsilon^{(i)}(w_2),\epsilon^{(i)}(w_1))$$
has only finitely many $\pm$'s (we read the signs from left to
right). For $w\in \mathscr{F}$, we define
\begin{equation*}
{\rm wt}(w)=\Lambda+\sum_{b\in\B}m_b\epsilon_b\in P,
\end{equation*}
where $m_b=\bigl|\{\,k\,|\,w_k=b\ (k\geq 1)
\,\}\bigr|-\delta_{|b|0}$ for $b\in \B$. Since $m_b=0$ for almost
all $b\in \B$, ${\rm wt}(w)$ is well-defined. For $w\in\mathscr{F}$
and $i\in I$, we set
\begin{equation*}
\varepsilon_i(w)={\rm max}\{\,k\,|\,e_i^kw\neq 0\,\}, \ \ \
\varphi_i(w)={\rm max}\{\,k\,|\,f_i^kw\neq 0\,\}.
\end{equation*}

For $c\in\mathbb{Z}$, let
\begin{equation*}
\begin{split}
H^c&=
\cdots \ \ov{c+3}\ \ov{c+2}\ \ov{c+1}, \ \ \ \ \ \  \text{if $c\geq 0$}, \\
&=\cdots \ov{3}\  \ov{2}\ \ov{1}\ \underbrace{1\ \cdots\ 1}_{|c|} \
, \ \ \ \ \ \ \ \ \ \ \text{if $c\leq 0$}.
\end{split}
\end{equation*}
Note that ${\rm wt}(H^0)=\Lambda$.
\begin{prop} $\mathscr{F}$
is a crystal graph for $\frak{g}$, and
\begin{equation*}
\mathscr{F}=\bigoplus_{c\in\mathbb{Z}}\B(c),
\end{equation*}
where $\B(c)$ is the connected component of $H^c$
$(c\in\mathbb{Z})$. Moreover,  $H^c$ is the highest weight element
in $\mathscr{F}$, that is, ${\rm wt}(H^c)\geq {\rm wt}(w)$ for
$w\in\B{(c)}$.
\end{prop}
\pf The conditions (b), (c) and (d) in Definition \ref{crystal
graph} are satisfied directly. So, it suffices to check that (a)
holds. Given $w=\cdots w_3 w_2 w_1 \in\mathscr{F}$ of charge $c$ and
$i\in I$, choose a sufficiently large $M>0$ such that
\begin{itemize}
\item[(i)] $w_{k}=\ov{c+k}$  for all $k\geq M$,

\item[(ii)] $\langle h_i, {\rm wt}(w^{> M})\rangle=0$, where $w^{> M}=\cdots w_{M+2}w_{M+1}\in\mathscr{F}$.
\end{itemize}
Put $w^{\leq M}=w_M w_{M-1}\cdots w_1$. Then we may view $w^{\leq
M}$ as an element in $\mathcal{W}_{m|n}$, where $m=M+c$ and $n\gg
0$. Also, we have $\varepsilon_i(w)=\varepsilon_i(w^{\leq M})$,
$\varphi_i(w)=\varphi_i(w^{\leq M})$, and $\langle h_i, {\rm
wt}(w^{\leq M})\rangle=\langle h_i, {\rm wt}(w)\rangle$. This
implies the condition (a).

Since $w^{\leq M}$ is $\frak{gl}_{m|n}$-equivalent to a semistandard
tableau $P(w^{\leq M})$ of a single column in $\B_{m|n}((1^M))$,
$w^{\leq M}$ is connected to the highest weight element
$H_{m|n}^{(1^M)}$. If $P(w^{\leq
M})=\rf_{i_1}\cdots\rf_{i_r}H_{m|n}^{(1^M)}$ for some $r\geq 0$ and
$i_k\in I_{m|n}$ ($1\leq k\leq r$), then we have
$w=\rf_{i_1}\cdots\rf_{i_r}H^c$. In particular, we have ${\rm
wt}(w)={\rm wt}(H^c)-\sum_{k=1}^r\alpha_{i_k}\leq {\rm wt}(H^c)$.
\qed\vskip 3mm

Let $x=\{\,x_b \,|\, b\in\B\, \}$. For $\mu\in P$, we define
$x^{\mu}=\prod_{b\in\B} x_b^{m_b}$, where
$\mu=k\Lambda+\sum_{b\in\B}m_b\epsilon_b$ ($k\in\mathbb{Z}$). Then
the character of $\mathscr{F}$ is given by
\begin{equation}
{\rm ch}\mathscr{F}=\sum_{w\in\mathscr{F}}x^{{\rm
wt}(w)}=\frac{\prod_{b\in\B^+}(1+x_b^{-1})}{\prod_{b'\in\B^-}(1-x_{b'})}.
\end{equation}

\subsection{Semi-infinite semistandard tableaux for $\frak{g}$}
Let us describe a crystal graph for $\frak{g}$ occurring as a
connected component in $\mathscr{F}^{\otimes u}$ ($u\geq 1$). Let
\begin{equation*}
\mathbb{Z}^u_+=\{\,\lambda=(\lambda_1,\cdots,\lambda_u)\in\mathbb{Z}^u\,|\,\lambda_1\geq
\lambda_2\geq\cdots\geq \lambda_u\,\}
\end{equation*}
be the set of all {\it generalized partitions} of length $u$.  For
each $\lambda=(\lambda_1,\cdots,\lambda_u)\in \mathbb{Z}^u_+$, we
call a $u$-tuple of semi-infinite words
$\w=(w^{(1)},\cdots,w^{(u)})\in\mathscr{F}^u$ a {\it semistandard
tableau of charge $\lambda$} if
\begin{itemize}
\item[(1)] $w^{(i)}=\cdots w^{(i)}_3w^{(i)}_2w^{(i)}_1\in\mathscr{F}$, where the charge of $w^{(i)}$ is $\lambda_i$ for $1\leq i\leq u$,

\item[(2)] $w^{(i)}_k\geq w^{(i+1)}_{k+d_i}$ for $1\leq i<u$ and $k\geq 1$,
where $d_i=\lambda_i-\lambda_{i+1}$,

\item[(3)] $w^{(i)}_k = w^{(i+1)}_{k+d_i}$ implies $|w^{(i)}_k|=0$.
\end{itemize}
We denote by $\B(\lambda)$ the set of all semistandard tableaux of
charge $\lambda$. In fact, each $\w\in\B(\lambda)$ determines a
unique semi-infinite tableau with infinitely many rows and $u$
columns, where each row of $\w$ reads (from left to right) as
follows;
\begin{equation*}
w^{(u)}_{k+d_1+\cdots+d_{u-1}}\cdots w^{(3)}_{k+d_1+d_2} w^{(2)}_{k+d_1}w^{(1)}_{k},
\end{equation*}
for $k\in\mathbb{Z}$ (we assume that $w^{(i)}_k$ is empty for $k\leq 0$).

\begin{ex}{\rm
Let $\lambda=(3,1,-2,-2)$, and let
$\w=(w^{(1)},w^{(2)},w^{(3)},w^{(4)})$ be given by

\begin{equation*}
\begin{array}{cccc}
 w^{(4)} & w^{(3)} & w^{(2)} & w^{(1)}  \\
 \vdots & \vdots &\vdots &\vdots   \\
 \ov{6} & \ov{6} & \ov{6} & \ov{6} \\
 \ov{5} & \ov{5} & \ov{5} & \ov{5} \\
 \ov{4} & \ov{4} & \ov{3} & \ov{1} \\
 \ov{3} & \ov{3} & \ov{2} &  \\
 \ov{2} & \ov{2} & \ov{1} &  \\
 \ov{1} & \ov{1} & &   \\
 1 & 2 & & \\
 1 & 3 & &
\end{array} \ \ .
\end{equation*}
Then $\w$ is a semistandard tableau of charge $\lambda$. }
\end{ex}

For $\lambda\in\mathbb{Z}^u_+$ and $\w=(w^{(1)},\cdots,
w^{(u)})\in\B(\lambda)$ , we may view
$\w=w^{(1)}\otimes\cdots\otimes w^{(u)}\in \mathscr{F}^{\otimes u}$,
and consider $\B(\lambda)$ as a subset of $\mathscr{F}^{\otimes u}$.

\begin{prop} For $\lambda\in\mathbb{Z}^u_+$,
$\B(\lambda)$ together with $0$ is stable under $e_i$ and $f_i$
$(i\in I)$. Hence, $\B(\lambda)$ is a crystal graph for $\frak{g}$.
Furthermore, $\B(\lambda)$ is a connected $I$-colored oriented graph
with a unique highest weight element $H^{\lambda}$.
\end{prop}
\pf To show that $\B(\lambda)$ is a crystal graph for $\frak{g}$, it
is enough to check that $\re_i\w,\rf_i\w \in \B(\lambda)\cup\{0\}$
for $\w\in\B(\lambda)$ and $i\in I$.

Suppose that $\w\in\B(\lambda)$ is given. Choose  sufficiently large
$m,n>0$. Then, for each $1\leq i\leq u$, we have
\begin{equation*}
w^{(i)}=\cdots \ov{m+3}\ \ov{m+2}\ \ov{m+1}\
w^{(i)}_*=H^{m}w^{(i)}_*,
\end{equation*}
for some $w^{(i)}_*\in\mathcal{W}_{m|n}$. Set
$\w_*=w^{(1)}_*\cdots w^{(u)}_*\in\mathcal{W}_{m|n}$. Then, $\w$
is $\frak{gl}_{m|n}$-equivalent to $\w_*$. Note that $P(\w_*)$,
the $P$-tableau of $\w_*$, is nothing but the $(m,n)$-hook
semistandard tableau which is obtained from $\w$ by removing the
entries smaller than $\ov{m}$. So, we have $\re_i\w,\rf_i\w \in
\B(\lambda)\cup\{0\}$ for $i\in I_{m|n}\subset I$.

Now, consider $H^{\lambda}=(w^{(1)},\cdots,w^{(u)})\in\B(\lambda)$,
where for $1\leq i\leq u$,
\begin{equation*}
\begin{split}
w^{(i)}& = H^{\lambda_i}, \hskip 4.3cm  \text{if $\lambda_i\geq 0$}, \\
       & = H^0\underbrace{(u-i+1)\cdots(u-i+1)}_{|\lambda_{i}|}, \ \ \  \ \ \text{if $\lambda_i< 0$}.
\end{split}
\end{equation*}
Since $P(\w_*)$ is connected to the highest weight element, it
follows that $\w$ is connected to $H^{\lambda}$.  Also, for
$\w\in\B(\lambda)$, we have ${\rm wt}(w)=\sum_{i=1}^r{\rm
wt}(w^{(i)})$ and ${\rm wt}(\w)\leq {\rm wt}(H^{\lambda})$. Hence,
$\B(\lambda)$ is connected with the unique highest weight element
$H^{\lambda}$. \qed

\begin{lem}\label{P1infty}
For $u\geq 1$, each connected component of $\mathscr{F}^{\otimes u}$
is isomorphic to $\B(\lambda)$ for some $\lambda\in\mathbb{Z}^u_+$.
\end{lem}

\pf Suppose that $\w=w^{(1)}\otimes\cdots\otimes
w^{(u)}\in\mathscr{F}^{\otimes u}$ is given. Choose  sufficiently
large $m,n>0$. Then, for each $1\leq i\leq u$, we have $
w^{(i)}=H^{m}w^{(i)}_*, $ for some $w^{(i)}_*\in\mathcal{W}_{m|n}$.

Set $\w_*=w^{(1)}_*\cdots w^{(u)}_*\in\mathcal{W}_{m|n}$. Then $\w$
is $\frak{gl}_{m|n}$-equivalent to $\w_*$, and
$P(\w_*)\in\B_{m|n}(\lambda)$ for some $\lambda=(\lambda_k)_{k\geq
1}\in\mathcal{P}_{m|n}$. We see that $\lambda_1\leq u$ from
Schensted's algorithm.

For $1\leq i\leq u$, let $w^{(i)}_{\sharp}$ be the word obtained by
reading the $(u-i+1)^{\rm th}$-column of $P(\w_*)$ from top to
bottom (note that the left-most column in $\lambda$ is the first
one). Set $\widetilde{\w}=\widetilde{w}^{(1)}\otimes\cdots\otimes
\widetilde{w}^{(u)}\in\mathscr{F}^{\otimes u}$ where
$\widetilde{w}^{(i)}=H^{m}w^{(i)}_{\sharp}$. Then  $\widetilde{\w}$
is $\frak{gl}_{m|n}$-equivalent to $\w$. Also, $\widetilde{\w}$ is
uniquely determined independent of all sufficiently large $m,n$, and
hence it is $\frak{g}$-equivalent to $\w$.

Let $\mu_i$ be the charge of $\widetilde{w}^{(i)}$ ($1\leq i\leq
u$). Since the length of $w_{\sharp}^{(i)}$ is less than or equal to
that of $w_{\sharp}^{(i+1)}$, we have $\mu_i\geq \mu_{i+1}$ for
$1\leq i\leq u-1$, and hence $\widetilde{\w}$ is a semistandard
tableau of charge $\mu$, where
$\mu=(\mu_1,\cdots,\mu_u)\in\mathbb{Z}^u_+$. \qed\vskip 3mm

By Lemma \ref{equivlence of tableaux}, we can also check the
following lemma, which implies that each $\w\in\mathscr{F}^{\otimes
u}$ is $\frak{g}$-equivalent to a unique semi-infinite semistandard
tableau.
\begin{lem}
Let $\w$ and $\w'$ be two semi-infinite semistandard tableaux. If
$\w\simeq_{\frak{g}}\w'$, then $\w=\w'$.\qed
\end{lem}

\subsection{Rational semistandard tableaux for $\frak{gl}_u$}
Let us recall the crystal graphs of rational representations of
$\frak{gl}_{u|0}$ for $u\geq 2$. By convention, we write
$\frak{gl}_u=\frak{gl}_{u|0}$, $\B_u=\B_{u|0}$, $I_u=I_{u|0}$,
$P_u=P_{u|0}$ and so on.

Let $\B_u^{\vee}= \{\, -\ov{1} < -\ov{2} < \cdots < -\ov{u} \}$ be
the dual crystal graph of $\B_u$ whose associated graph is given by
\begin{equation*}
-\ov{1}\ \stackrel{\ov{1}}{\longrightarrow}\ -\ov{2}\
\stackrel{\ov{2}}{\longrightarrow}
\cdots\stackrel{\ov{u-1}}{\longrightarrow}\ -\ov{u},
\end{equation*}
where ${\rm wt}(-\ov{k})=-{\rm wt}(\ov{k})=-\epsilon_{\ov{k}}$ for
$1\leq k\leq u$ (cf.\cite{Kas94}).

Given $\lambda=(\lambda_1,\cdots,\lambda_u)\in \mathbb{Z}_+^u$, we
may identify $\lambda$ with a {\it generalized Young diagram} in the
following way. First, we fix a vertical line. Then for each
$\lambda_{k}$, we place $|\lambda_k|$ nodes (or boxes) in the
$k^{\rm th}$ row in a left-justified (resp. right-justified) way
with respect to the vertical line if $\lambda_{k}\geq 0$ (resp.
$\lambda_k\leq 0$). For example,\vskip 3mm
\begin{center}
 $\lambda=(3,2,0,-1,-2) \longleftrightarrow$
\begin{tabular}{cc|ccc}
   &   & $\bullet$ & $\bullet$ & $\bullet$  \\
   &   & $\bullet$ & $\bullet$ &   \\
   & & & &  \\
   & $\bullet$ &   &   &  \\
 $\bullet$ & $\bullet$ &   & &   \\
\multicolumn{5}{c}{\small -2  \ -1 \ \  1 \ \ 2 \ \ 3  }
\end{tabular}\ \ \ \ .
\end{center}
We enumerate the columns of a diagram as in the above figure.

\begin{df}[cf.\cite{St}]\label{rationalSST}{\rm
Let  $T$ be a tableau obtained by filling a generalized Young
diagram $\lambda$ of length $u$ with the entries in $\B_u \cup
\B_u^{\vee}$. We call $T$  a {\it rational semistandard of shape
$\lambda$} if
\begin{itemize}
\item[(1)] the entries in the columns indexed by positive (resp. negative)
numbers belong to $\B_u$ (resp. $\B_u^{\vee}$),
\item[(2)] the entries in each row (resp. column) are weakly (resp. strictly)
increasing from left to right (resp. from top to bottom),
\item[(3)] if $b_1<\cdots<b_s$ (resp. $-b'_1<\cdots<-b'_t$) are the entries
in the $1^{\rm st}$ (resp. $-1^{\rm st}$) column ($s+t\leq u$), then
$$b''_i\leq b_i,$$ for $1\leq i\leq s$, where $\{ b''_1<\cdots<b''_{u-t} \}=
\B_u\setminus \{ b'_1,\cdots,b'_t \}$.
\end{itemize}}
\end{df}
We denote by $\B_u(\lambda)$ the set of all rational semistandard
tableaux of shape $\lambda$.

\begin{ex}{\rm For $\lambda=(3,2,0,-1,-2)$, we have

\begin{center}
\begin{tabular}{cc|ccc}
   &   & $\ov{4}$ & $\ov{3}$ & $\ov{1}$ \\
   &   & $\ov{2}$ & $\ov{2}$ &  \\
   & & & & \\
   & $-\ov{1}$ &   &   &  \\
 $-\ov{2}$ & $-\ov{2}$ &   &   &
\end{tabular}$\in \B_5(\lambda)$,\ \  but
\begin{tabular}{cc|ccc}
   &   & $\ov{4}$ & $\ov{3}$ & $\ov{1}$ \\
   &   & $\ov{2}$ & $\ov{2}$ &  \\
   & & & & \\
   & $-\ov{4}$ &   &   &  \\
 $-\ov{2}$ & $-\ov{5}$ &   &   &
\end{tabular}$\not\in \B_5(\lambda).$
\end{center}}
\end{ex}

Let us explain the relation between the crystal graphs of rational
representations and polynomial representations. Let $T$ be a
rational semistandard tableau in $ \B_u((-1)^t)$ ($0\leq t\leq u$)
with the entries $-b_1<\cdots<-b_t$. We define $\sigma(T)$ be the
tableau in $\B_u(1^{u-t})$ with the entries $b'_1<\cdots<b'_{u-s}$,
where $\{ b'_1<\cdots<b'_{u-t}\}= \B_u\setminus \{b_1<\cdots<b_t\}$.
If $t=u$, then we define $\sigma(T)$ to be the empty word.

Generally, for $\lambda\in\mathbb{Z}_+^u$ and $T\in\B_u(\lambda)$,
we define $\sigma(T)$ to be the tableau obtained by applying
$\sigma$ to the $-1^{\rm st}$ column of $T$. For example, when
$u=5$, we have \vskip 3mm
\begin{center}
$\sigma\left(\begin{array}{cc|ccc}
   &   & \ov{4} & \ov{3} & \ov{1} \\
   &   & \ov{2} & \ov{2} &  \\
   & & & & \\
   & {\bf -\ov{1}} &   &   &  \\
 -\ov{2} & {\bf -\ov{4}} &   &   &
\end{array}\right)$\ \ =
\begin{tabular}{c|cccc}
   & {\bf $\bf \ov{5}$} & $\ov{4}$ & $\ov{3}$ & $\ov{1}$ \\
   & {\bf $\bf \ov{3}$} & $\ov{2}$ & $\ov{2}$ &  \\
   & {\bf $\bf \ov{2}$} & & & \\
   &  &   &   &  \\
 $-\ov{2}$ &  &   &   &
\end{tabular}.
\end{center}\vskip 3mm
By Definition \ref{rationalSST}, we have
$\sigma(T)\in\B_u(\lambda+(1^u))$, where we define
$\mu+\nu=(\mu_k+\nu_k)_{k\geq 1}$ for two generalized partitions
$\mu=(\mu_k)_{k\geq 1}$ and $\nu=(\nu_k)_{\geq 1}$ in
$\mathbb{Z}^u_+$.

Now, we embed $\B_u(\lambda)$ into $(\B_u\oplus\B_u^{\vee})^{\otimes
N}$ ($N=\sum_{k\geq 1}|\lambda_k|$) by column reading of tableaux
(cf. Section 2.2), and apply $e_i,f_i$ to $\B_u(\lambda)$ ($i\in
I_u$). Then
\begin{lem}\label{shift of SST}
For $\lambda\in\mathbb{Z}^u_+$, $\B_u(\lambda)$ is a crystal graph
for $\frak{gl}_u$, and the map
$$\sigma : \B_u(\lambda)\cup\{0\} \rightarrow \B_u(\lambda+(1^u))\cup\{0\},$$
$(\sigma(0)=0)$ is a bijection which commutes with $e_i,f_i$ for
$i\in I_u$, where ${\rm wt}(\sigma(T))={\rm
wt}(T)+(\epsilon_{\ov{u}}+\cdots+\epsilon_{\ov{1}})$ for $T\in
\B_u(\lambda)$.
\end{lem}
\pf It follows immediately that $\sigma$ is a bijection and ${\rm
wt}(\sigma(T))={\rm
wt}(T)+(\epsilon_{\ov{u}}+\cdots+\epsilon_{\ov{1}})$ for $T\in
\B_u(\lambda)$.

So, it remains to show that $\B_u(\lambda)$ is a crystal graph for
$\frak{gl}_u$, and $\sigma$ commutes with the Kashiwara operators.
We will use induction on  the number of columns indexed by negative
numbers.

Suppose that $\lambda=((-1)^t)$ for some $0\leq t\leq u$. Then it is
straightforward to check that  our claim holds.

Next, for a general $\lambda\in\mathbb{Z}^u_+$, any tableau  $T$ in
$\B_u(\lambda)$ can be viewed as a tensor product of its columns
when we apply $e_i,f_i$. By Definition \ref{tensor product} and the
argument in the case of a single column, it follows that
$\sigma(x_iT)=x_i\sigma(T)$ for $x=e,f$ and $i\in I_u$. Also, by
induction hypothesis, $x_i\sigma(T)\in \B_u(\lambda+
(1^u))\cup\{0\}$, which is equivalent to saying that $x_iT\in
\B_u(\lambda)\cup\{0\}$. This completes the induction. \qed\vskip
3mm

Note that for $\lambda\in\mathbb{Z}_+^u$, there exists a unique
highest weight element $H^{\lambda}_u$ in $\B_u(\lambda)$. In fact,
$H^{\lambda}_u=\sigma^{-k}(H_u^{\lambda+(k^u)})$ for all $k\geq 0$
such that $\lambda+(k^u)$ is an ordinary partition.

Let $y=\{\,y_b\,|\, b\in \B_u\,\}$ be the set of variables indexed
by $\B_u$. For $\mu=\sum_{b\in \B_u}\mu_b\epsilon_b\in P_u$, we set
$y^{\mu}=\prod_{b\in \B_u}y_b^{\mu_b}$. For
$\lambda\in\mathbb{Z}_+^u$, the character of $\B_u(\lambda)$ is
given by a {\it rational Schur function} corresponding to $\lambda$;
\begin{equation*}
s_{\lambda}(y)=\sum_{T\in\B_u(\lambda)}y^{{\rm wt}T}.
\end{equation*}
By Lemma \ref{shift of SST}, we have
$s_{\lambda+(1^u)}(y)=(y_{\ov{u}}\cdots
y_{\ov{1}})s_{\lambda}(y)$.

\subsection{Decomposition of $\mathscr{F}^{\otimes u}$}

Now, let us decompose $\mathscr{F}^{\otimes u}$ for $u\geq 2$. We
start with another description of $\mathscr{F}^{\otimes u}$ in terms
of matrices of non-negative integers. Set
\begin{equation}
\begin{split}
\mathscr{M}^u=\{ &A=(a_{bb'})_{b\in\B, b'\in \B_{u}}\,|\, \\
& \text{(1) $a_{bb'}\in\mathbb{Z}_{\geq 0}$,} \\
& \text{(2) $a_{bb'}\leq 1$ if $|b|=0$,} \\
& \text{(3) $a_{bb'}=1$ for all $b\ll \ov{1}$, and $a_{bb'}=0$ for
all $b\gg 1$}\, \}.
\end{split}
\end{equation}

For $\w=w^{(1)}\otimes\cdots\otimes w^{(u)}\in\mathscr{F}^{\otimes
u}$, set $A(\w)=(a_{b\, \overline{k}})\in\mathscr{M}^u$, where
$a_{b\, \overline{k}}$ is the number of occurrences of $b$ in
$w^{(u-k+1)}$ ($1\leq k\leq u$). Then the map $\w\mapsto A(\w)$ is a
bijection from $\mathscr{F}^{\otimes u}$ to $\mathscr{M}^u$, where
each $w^{(i)}$ corresponds to the ${i}^{\rm th}$-column of $A(\w)$
for $1\leq i\leq u$.\vskip 3mm

\begin{ex}\label{correspondence}{\rm
\begin{equation*}
\mathscr{F}^{\otimes 3}\ni \w=\begin{array}{c}
 \vdots \\
 \ov{5} \\
 \ov{4} \\
 \ov{2} \\
 2  \\
 2 \\
 \mbox{} \\
 \mbox{}
\end{array}
\otimes
\begin{array}{c}
 \vdots \\
 \ov{5} \\
 \ov{4} \\
 \ov{3} \\
 3 \\
  \mbox{} \\
  \mbox{} \\
  \mbox{}
\end{array}
\otimes
\begin{array}{c}
 \vdots \\
 \ov{5} \\
 \ov{3} \\
 \ov{2} \\
  1 \\
  1 \\
  1 \\
  \mbox{}
\end{array}
\ \ \ \longleftrightarrow\ \ \ A(\w)=\left(\begin{tabular}{ccc}
 \vdots & \vdots & \vdots \\
 1 & 1 & 1  \\
 1 & 1 & 0  \\
 0 & 1 & 1  \\
 1 & 0 & 1  \\
 0 & 0 & 0  \\ \hline
 0 & 0 & 3  \\
 2 & 0 & 0  \\
 0 & 1 & 0  \\
 \vdots & \vdots & \vdots
\end{tabular}\right )\in\mathscr{M}^3.
\end{equation*}
}
\end{ex}

Let us define the crystal graph structures on $\mathscr{M}^u$, which
are naturally induced from those on $\M^{\sharp}_{m|n,u|0}$ in
Section 4. For $m,n>0$, we define
$${\rm Res}_{m|n} : \mathscr{M}^u \longrightarrow \M^{\sharp}_{m|n,u|0},$$
by ${\rm Res}_{m|n}A=(a_{bb'})_{b\in\B_{m|n},b'\in\B_{u}}$ for
$A=(a_{bb'})\in  \mathscr{M}^u$.\vskip 3mm

For $A\in \mathscr{M}^u$ and $i\in I$, we define $e_iA$ and $f_iA$
to be the unique elements in $\mathscr{M}^u \cup \{0\}$ satisfying
\begin{equation*}
{\rm Res}_{m|n}(e_iA)=e_i({\rm Res}_{m|n}A), \ \ {\rm
Res}_{m|n}(f_iA)=f_i({\rm Res}_{m|n}A),
\end{equation*}
for all sufficiently large $m,n>0$, where we assume that ${\rm
Res}_{m|n}0=0$. We set
\begin{equation*}
\begin{split}
&{\rm wt}(A)=u\Lambda+\sum_{b\in\B}m_b\epsilon_b,\\
&\varepsilon_i(A)={\rm max}\{\,k\,|\,e_i^kA\neq 0\,\}, \ \ \
\varphi_i(A)={\rm max}\{\,k\,|\,f_i^kA\neq 0\,\},
\end{split}
\end{equation*}
for $A\in\mathscr{M}^u$ and $i\in I$, where
$m_b=\sum_{b'\in\B_u}(a_{bb'}-\delta_{|b|0})$. Then it is not
difficult to check the following lemma;
\begin{lem}\label{P1infty'}
$\mathscr{M}^u$ is a crystal graph for $\frak{g}$, and the map
$\w\mapsto A(\w)$ is an isomorphism of $\frak{g}$-crystals from
$\mathscr{F}^{\otimes u}$ to $\mathscr{M}^u$. In particular, for
$A\in\mathscr{M}^u$, $A$ is $\frak{g}$-equivalent to a unique
semistandard tableau in $\B(\lambda)$ for some
$\lambda\in\mathbb{Z}^u_+$. \qed
\end{lem}\vskip 3mm

Next, let us define a $\frak{gl}_u$-crystal structure on
$\mathscr{M}^u$. Since $\mathscr{M}^u$ can be identified with
$\mathscr{F}^{\otimes u}$, this induces a $\frak{gl}_u$-crystal
structure on $\mathscr{F}^{\otimes u}$. For $A\in \mathscr{M}^u$ and
$j\in I_{u}$, we define $e_j^*A$ and $f_j^*A$ to be the unique
elements in $\mathscr{M}^u \cup \{0\}$ satisfying
\begin{equation*}
{\rm Res}_{m|n}(e_j^*A)=e_j^*({\rm Res}_{m|n}A), \ \ {\rm
Res}_{m|n}(f_j^*A)=f_j^*({\rm Res}_{m|n}A),
\end{equation*}
for all sufficiently large $m,n>0$. We set
\begin{equation*}
\begin{split}
&{\rm wt}^*(A)= \sum_{b\in\B_u}m_b\epsilon_b, \\
&\varepsilon^*_j(A)={\rm max}\{\,k\,|\,(e^*_j)^kA\neq 0\,\}, \ \ \
\varphi^*_j(A)={\rm max}\{\,k\,|\,(f^*_j)^kA\neq 0\,\},
\end{split}
\end{equation*}
for $A\in\mathscr{M}^u$ and $j\in I_u$, where
$m_b=\sum_{b'\in\B}(a_{b'b}-\delta_{|b'|0})$. Then
\begin{lem}\label{P2infty}
$\mathscr{M}^u$ is a crystal graph for $\frak{gl}_u$ with respect to
$e_j^*,f_j^*$ $(j\in I_u)$, and each $A\in\mathscr{M}^u$ is
$\frak{gl}_u$-equivalent to a unique rational semistandard tableau.
\end{lem}
\pf It is clear that $\mathscr{M}^u$ is a crystal graph for
$\frak{gl}_u$. For $A\in\mathscr{M}^u$, choose sufficiently large
$m,n$. As an element in a $\frak{gl}_u$-crystal, consider a
semistandard tableau $T=P_2^{\sharp}({\rm Res}_{m|n}A)$  (see
\eqref{pisharp}). We observe that ${\rm wt}(T)={\rm
wt}^*(A)+k(\epsilon_{\ov{u}}+\cdots+\epsilon_{\ov{1}})$ for some
$k\geq 0$. Then, by Lemma \ref{shift of SST}, the rational
semistandard tableau $\sigma^{-k}(T)$, say $S$, is
$\frak{gl}_u$-equivalent to $A$. Note that $S$ is independent of $m$
and $n$ since ${\rm wt}^*(A)$ is fixed. \qed\vskip 3mm

Suppose that $A\in\mathscr{M}^u$ is given. By Lemma \ref{P1infty'},
there exists a unique semistandard tableau $\w$ in $\B(\lambda)$ for
some $\lambda\in\mathbb{Z}_+^u$, which is $\frak{g}$-equivalent to
$A$, and we set $\mathscr{P}_1(A)=\w$. Also, by Lemma \ref{P2infty},
there exists a unique rational semistandard tableau $T$, which is
$\frak{gl}_u$-equivalent to $A$, and we set $\mathscr{P}_2(A)=T$.
Now, we define
\begin{equation}\label{piinf}
\varpi(A)=(\mathscr{P}_1(A),\mathscr{P}_2(A)).
\end{equation}

\begin{prop}\label{bicrystal''}
$\mathscr{M}^u$ is a $(\frak{g},\frak{gl}_{u})$-bicrystal, and for
each connected component $\mathscr{C}$ in $\mathscr{M}^u$, $\varpi$
gives the following isomorphism of
$(\frak{g},\frak{gl}_{u})$-bicrystals;
$$\varpi : \mathscr{C} \longrightarrow \B(\lambda)\times \B_u(\mu),$$
for some $\lambda,\mu\in\mathbb{Z}_+^u$.
\end{prop}
\pf It follows from Proposition \ref{bicrystal'}. \qed\vskip 3mm

\begin{ex}{\rm
Let $A$ be given in Example \ref{correspondence}. Then
\begin{equation*}
A\ \simeq_{\frak{g}} \
\begin{array}{c}
 \vdots \\
 \ov{5} \\
 \ov{4} \\
 \ov{2} \\
 2  \\
 2 \\
 \mbox{} \\
 \mbox{}
\end{array}
\otimes
\begin{array}{c}
 \vdots \\
 \ov{5} \\
 \ov{4} \\
 \ov{3} \\
 3 \\
  \mbox{} \\
  \mbox{} \\
  \mbox{}
\end{array}
\otimes
\begin{array}{c}
 \vdots \\
 \ov{5} \\
 \ov{3} \\
 \ov{2} \\
  1 \\
  1 \\
  1 \\
  \mbox{}
\end{array}\
\simeq_{\frak{g}}\
\begin{array}{ccc}
 \vdots & \vdots & \vdots \\
 \ov{5} & \ov{5} & \ov{5} \\
 \ov{4} & \ov{4} & \ov{2} \\
 \ov{3} & \ov{3}   \\
\ov{2} & 2 \\
 1 & 2 \\
 1 & 3 \\
 1
\end{array}\ =\ \mathscr{P}_1(A)\in\B(3,-1,-2),
\end{equation*}
On the other hand, consider ${\rm Res}_{4|3}A$. Then
\begin{equation*}
\begin{split}
{\rm Res}_{4|3}A\ &\simeq_{\frak{gl}_u} \ \ov{2}\ \otimes \ov{3}
\, \ov{3} \otimes \ov{1} \, \ov{1} \, \ov{1} \otimes
\begin{array}{c}
 \ov{3} \\
 \ov{1}
\end{array}\otimes
\begin{array}{c}
 \ov{2} \\
 \ov{1}
\end{array}\otimes
\begin{array}{c}
 \ov{3} \\
 \ov{2}
\end{array} \\
& \simeq_{\frak{gl}_u}\
\begin{array}{ccccccc}
\ov{3} & \ov{3} &\ov{3} &\ov{3} &\ov{2} &\ov{1} & \\
\ov{2} & \ov{2} &\ov{1} &\ov{1} & \ov{1} & & \\
\ov{1} & & & & & &
\end{array}\ = \ P_2^{\sharp}({\rm Res}_{4|3}A).
\end{split}
\end{equation*}
Since ${\rm wt}(A)=-\epsilon_{\ov{2}}+\epsilon_{\ov{1}}$, we have
\begin{equation*}
\mathscr{P}_2(A)=\sigma^{-4}(\ P_2^{\sharp}({\rm Res}_{4|3}A))\ = \
\begin{array}{ccc|ccc}
 & & & \ov{2} &\ov{1} &\\
 & & & \ov{1} &  & \\
 -\ov{1} & -\ov{2} & -\ov{2} & & &
\end{array}\in \B_3(2,1,-3).
\end{equation*} }
\end{ex}\vskip 3mm

Now, we let
\begin{equation}
\mathscr{M}^u_{\rm h.w.}=\{\,A\,|\,\varpi(A)=(H^{\lambda},H^{\mu}_u)
\ \ \text{for some $\lambda,\mu\in\mathbb{Z}^u_+$ }\,\}
\end{equation}
be the set of all the highest weight elements in $\mathscr{M}^u$.
Then $\mathscr{M}^u$ is the direct sum of the connected components
of the elements in $\mathscr{M}^u_{\rm h.w.}$ as a
$(\frak{g},\frak{gl}_{u})$-bicrystal. For
$\lambda\in\mathbb{Z}^u_+$, let $\mu_i$ ($1\leq i\leq u$) be the
number of occurrences of $i$ in $H^{\lambda}$. Put
$\nu_i=\mu_i-\mu_{i+1}$ for $1\leq i\leq u$ where $\mu_{u+1}=0$.
Then we define ${\mathscr A}_{\lambda}=(a_{bb'})\in\mathscr{M}^u$ by
\begin{itemize}
\item[(1)] for $k\geq 1$, $0\leq l< u$,
\begin{equation*}
a_{\ov{k}\, \ov{u-l}}=
\begin{cases}
1, & \text{if there exists $\ov{k}$ in the $(l+1)^{\rm th}$-column of $H^{\lambda}$}, \\
0, & \text{otherwise}.
\end{cases}
\end{equation*}

\item[(2)] for $k\geq 1$, $0\leq l< u$,
\begin{equation*}
a_{{k}\, \ov{u-l}}=
\begin{cases}
\nu_{k+l}, & \text{if $1\leq k+l \leq u$}, \\
0, & \text{otherwise},
\end{cases}
\end{equation*}
\end{itemize}
(cf.\eqref{Alambda2}). Then, we can check that $\varpi({\mathscr
A}_{\lambda})=(H^{\lambda},H^{\lambda^*}_u)$, where
$$\lambda^*=(-\lambda_u,\cdots,-\lambda_1)\in\mathbb{Z}^u_+.$$
\begin{thm}\label{highest weight vector3} We have
\begin{equation*}
\mathscr{M}^u_{\rm
h.w.}=\{\,\mathscr{A}_{\lambda}\,|\,\lambda\in\mathbb{Z}^u_+\,\},
\end{equation*}
and the following isomorphism of
$(\frak{g},\frak{gl}_{u})$-bicrystals;
\begin{equation*}
\varpi : \mathscr{M}^{u} \longrightarrow
\bigoplus_{\lambda\in\mathbb{Z}^u_+} \B(\lambda)\times
\B_{u}(\lambda^*).
\end{equation*}
\end{thm}\pf Let $A$ be a highest weight element in $\mathscr{M}^u$.
Then, ${\rm Res}_{m|n}A$ is a highest weight element in
$\M^{\sharp}_{m|n,u|0}$ for all sufficiently large $m,n>0$. By
Theorem \ref{highest weight vector2}, it follows that
$A=\mathscr{A}_{\lambda}$ for some $\lambda\in\mathbb{Z}^u_+$.
\qed\vskip 3mm

The character of $\mathscr{M}^u$ is given by
\begin{equation*}
{\rm ch}\mathscr{M}^u =\sum_{A\in\mathscr{M}^u}x^{{\rm
wt}(A)}y^{{\rm
wt}^*(A)}=\frac{\prod_{b\in\B^+}\prod_{b'\in\B_{u}}(1+x_{b}^{-1}y_{b'}^{-1}
)} {\prod_{b\in\B^-}\prod_{b'\in\B_{u}}(1-x_{b}y_{b'})}.
\end{equation*}
By Theorem \ref{highest weight vector3}, we obtain the following
identity;
\begin{equation}\label{infiniteCauchy}
\frac{\prod_{b\in\B^+}\prod_{b'\in\B_{u}}(1+x_{b}^{-1}y_{b'}^{-1}
)} {\prod_{b\in\B^-}\prod_{b'\in\B_{u}}(1-x_{b}y_{b'})}
=\sum_{\lambda\in\mathbb{Z}^u_+}{\rm
ch}\B(\lambda)s_{\lambda^*}(y),
\end{equation}
where ${\rm ch}\B(\lambda)=\sum_{\w\in\B(\lambda)}x^{{\rm wt}(\w)}$
is the character of $\B(\lambda)$, and $s_{\lambda^*}(y)$ is the
rational Schur function corresponding to $\lambda^*$.

From the classical Cauchy identities of Schur functions
(cf.\cite{Mac}), it follows that the left hand side of
\eqref{infiniteCauchy} is equal to
\begin{equation*}
\sum_{\mu,\nu\in\mathcal{P}}s_{\mu}(y^{-1})s_{\mu'}(x_+^{-1})
s_{\nu}(y)s_{\nu }(x_-),
\end{equation*}
where $\mathcal{P}$ is the set of all partitions,
$y^{-1}=\{\,y_b^{-1}\,|\,b\in\B_u\,\}$,
$x_{+}^{-1}=\{\,x_b^{-1}\,|\,b\in\B^+\,\}$, and
$x_{-}=\{\,x_b\,|\,b\in\B^-\,\}$. Note that
$s_{\mu}(y^{-1})=s_{\mu^*}(y)$ and
\begin{equation*}
s_{\mu^*}(y)s_{\nu}(y)=\sum_{\lambda\in\mathbb{Z}^u_+}N^{\lambda}_{\mu\nu}s_{\lambda}(y),
\end{equation*}
for some $N^{\lambda}_{\mu\nu}\in\mathbb{Z}_{\geq 0}$. Then the
left hand side of \eqref{infiniteCauchy} can be written as
\begin{equation*}
\sum_{\lambda\in\mathbb{Z}^u_+}\left(\sum_{\mu,\nu\in\mathcal{P}}
N^{\lambda}_{\mu\nu}s_{\mu'}(x_+^{-1})s_{\nu}(x_-)\right)s_{\lambda}(y).
\end{equation*}
From a linear independence of
$\{\,s_{\nu}(y)\,|\,\nu\in\mathbb{Z}^u_+\,\}$ (see Lemma 3.1
\cite{CL}), we obtain the following character formula of
$\B(\lambda)$  by comparing with the right hand side of
\eqref{infiniteCauchy};
\begin{cor}[cf.\cite{CL,KacR}] For $\lambda\in\mathbb{Z}^u_+$, we
have
\begin{equation*}
{\rm ch}\B(\lambda)=\sum_{\mu,\nu\in\mathcal{P}} N^{\lambda^*}_{\mu
\nu}s_{\mu'}(x_+^{-1})s_{\nu}(x_-).
\end{equation*}
\end{cor}\qed

\end{document}